\documentclass{amsart}

\theoremstyle{plain}
\newtheorem{theorem}{Theorem}[section]
\newtheorem{proposition}[theorem]{Proposition}
\newtheorem{corollary}[theorem]{Corollary}
\newtheorem{lemma}[theorem]{Lemma}

\newtheorem{definition}[theorem]{Definition}

\newtheorem{remark}[theorem]{Remark}

\newcommand{\bfu}{{\bf u}}

\newcommand{\bfC}{{\mathbb C}}

\newcommand{\bfR}{{\mathbb R}}

\newcommand{\bari}{{\overline i}}
\newcommand{\barj}{{\overline j}}

\newcommand{\barpartial}{{\overline \partial}}

\newcommand{\barp}{{\overline p}}
\newcommand{\barq}{{\overline q}}
\newcommand{\barz}{{\overline z}}
\newcommand{\barv}{{\overline v}}
\newcommand{\barw}{{\overline w}}

\newcommand{\tildeX}{\widetilde {X}}

\newcommand{\tildeg}{{\widetilde g}}

\newcommand{\tildex}{{\widetilde x}}

\newcommand{\tildeomega}{{\widetilde \omega}}

\newcommand{\mapright}[1]{\smash{\mathop{   \hbox to 0.7cm{\rightarrowfill}}
  \limits^{#1}}}

\newcommand{\Ker}{{\rm Ker}}

\def \R {\mathbb R}
\def \cal {\mathcal}
\def \vs{\vspace*{0.1cm}} 

\def \ds{\displaystyle}
\def\SS{S}
\def\DD{D}
\def\p{\partial}
\def\wt{\widetilde}

\def\D{\Delta}
\def\p{\partial}
\def\l{\lambda}

\def \ds{\displaystyle}
\def\a{\alpha}
\def\b{\beta}
\def\t{\theta}

\def\cal{\mathcal}

\def\n{\nabla}
\def\O{\Omega}

\def\C{\mathbb C}
\def\la{\langle}
\def\ra{\rangle}
\def\l{\lambda}
\def\L{\Lambda}

\def\part{\partial}
\def\R{{\mathbb R}}

\def\g{{\mathfrak g}} 
 
\def\hh{\mathfrak h}
\def\o{\omega}
\def\ba{\begin{array}}
\def\ea{\end{array}}

\newtheorem{prop}[theorem]{Proposition}

\newtheorem{defn}[theorem]{Definition}

\newtheorem{rem}[theorem]{Remark}
\newtheorem{lem}[theorem]{Lemma}

\newtheorem{pro}[theorem]{Proposition}

\def\vol{\mathop{\mathrm{Vol}}\nolimits}

\def\aut{\mathop{\mathrm{Aut}}\nolimits}

\def\aut{\mathop{\mathrm{Aut}}\nolimits}

\def\KCM{\mathop{\mathrm{KCM}}\nolimits}
\def\isom{\mathop{\mathrm{Isometry}}\nolimits}

\begin{document}

\title
{Transverse K\"ahler geometry of Sasaki manifolds and toric Sasaki-Einstein manifolds}
\author{Akito Futaki}
\address{Department of Mathematics, Tokyo Institute of Technology, 2-12-1,
O-okayama, Meguro, Tokyo 152-8551, Japan}
\email{futaki@math.titech.ac.jp}
\author{Hajime Ono}
\address{Department of Mathematics, Tokyo Institute of Technology, 2-12-1,
O-okayama, Meguro, Tokyo 152-8551, Japan}
\email{ono@math.titech.ac.jp}
\author{Guofang Wang}
\address{Max Planck Institute for Mathematics in the Sciences
Inselstr. 22-26, 04103 Leipzig, Germany
}
\email{gwang@mis.mpg.de}

\date{December 4, 2006 }

\begin{abstract} 
In this paper we study compact Sasaki manifolds in view of transverse K\"ahler geometry and
extend some results in K\"ahler geometry to Sasaki manifolds. In particular we 
define integral invariants which obstruct the existence of transverse K\"ahler metric with
harmonic Chern forms.  The integral invariant $f_1$ for the first Chern class case becomes an obstruction to the 
existence of transverse K\"ahler metric of constant scalar curvature. 
We prove the existence of transverse 
K\"ahler-Ricci solitons (or {\it Sasaki-Ricci soliton}) on
compact toric Sasaki manifolds whose basic first Chern form of the normal bundle of the Reeb
foliation is positive and the first Chern class of the contact bundle is trivial.
We will further show that if $S$ is a compact toric Sasaki manifold with the above assumption
then by deforming the Reeb field we get a 
Sasaki-Einstein structure on $S$.
As an application we obtain irregular toric Sasaki-Einstein metrics on 
the unit circle bundles
of the powers of the canonical bundle of the two-point blow-up of the complex projective plane.
\end{abstract}
\keywords{Sasaki manifold,  Einstein metric, transverse K\"ahler geometry}

\subjclass{Primary 53C55, Secondary 53C21, 55N91 }

\maketitle

\section{Introduction}
Sasaki manifolds can be studied from many view points as they have many
structures. They are characterized by having a K\"ahler structure on the 
cone. They have a one dimensional foliation, called the Reeb foliation,
which has a transverse K\"ahler structure. They also have a contact structure,
which provides us a moment map. These structures are intimately related each other, but when we
study the deformations of Sasaki structures it is efficient to fix some of the structures
and vary other structures. In this paper we study the deformations of Sasaki structure
fixing the Reeb foliation together with its transverse holomorphic structure and 
the holomorphic structure of the cone, while varying the K\"ahler metric on the
transverse holomorphic structure and, as a result, the contact structure. 

Sasaki geometry is often described as an odd dimensional analogue of K\"ahler
geometry. The above deformations on a Sasaki manifold correspond to the deformations of 
K\"ahler forms in a fixed K\"ahler class on a K\"ahler manifold. Therefore we may
try to extend results related to Calabi's extremal problem in K\"ahler geometry 
to the above setting in Sasaki geometry. The normal bundle of the Reeb foliation has basic
Chern forms which are expressed by basic differential forms. The basic forms are those differential
forms obtained by pulling back differential forms on the local leaf spaces of the Reeb foliation, 
see section 4 for the precise definition.
The basic first Chern class $c_1^B$ of the normal bundle of the Reeb foliation is said to 
be positive if it is represented by a positive basic 2-form
in the sense of K\"ahler geometry on the local leaf spaces; 
We will describe this condition by the notation $c_1^B > 0$. 
When one considers the problem of finding a Sasaki-Einstein metric, which is one of
the main interests of this paper,  it is reduced to
the problem of finding a transverse K\"ahler-Einstein metric.
A Sasaki-Eisntein metric must satisfy $Ric = 2m\,g$ if $\dim S = 2m+1$, and then
the transverse K\"ahler-Einstein metric satisfies
$$\rho^T = (2m+2) \omega^T$$
where $\omega^T$ and $\rho^T$ are respectively the transverse K\"ahler form and
the transverse Ricci form. See Section 3 below.
Since $\rho^T$ represents the basic first Chern class $c_1^B$ and 
$\omega^T$ is given by $\frac12 d\eta$, $\eta$ being the contact $1$-form dual to the
Reeb vector field, one is naturally lead to another assumption
$c_1(D) = 0$ as a de Rham cohomology class where $D = \mathrm{Ker}\, \eta$ is the contact bundle, see Proposition \ref{c_1(D)}
 for the detail.
These two conditions $c_1^B > 0$ and $c_1(D) = 0$ 
are primary obstructions for the existence of a transverse K\"ahler-Einstein metric
of positive scalar curvature, or equivalently for the existence of a Sasaki-Einstein metric.
Proposition \ref{c_1(D)} asserts that the two conditions $c_1^B > 0$ and $c_1(D) = 0$ are 
equivalent to that
$c_1^B$ is represented by $\tau d\eta$ for some positive constant $\tau$.
By changing $\eta$ by homothety we may then assume that $c_1^B = (2m+2) [\omega^T]$.

In this paper we first extend obstructions to the existence of K\"ahler
metric  of harmonic Chern forms (\cite{bando83}) to the transverse K\"ahler geometry
of compact Sasaki manifolds. 
The invariant $f_1$ for the first Chern form is an obstruction to the existence
of transverse K\"ahler metric of constant scalar curvature, which is a
secondary obstruction to the
existence of transverse K\"ahler-Einstein metric when $c_1^B > 0$ and $c_1(D) = 0$. This extension has been obtained recently by Boyer, Galicki and Simanca (\cite{BGS})
independently.

More generally than transverse K\"ahler-Einstein metrics of positive scalar curvature 
we consider the transverse K\"ahler-Ricci solitons, or
{\it Sasaki-Ricci solitons} for short. A Sasaki metric is called a Sasaki-Ricci
soliton if there exists a Hamiltonian holomorphic vector field $X$ such that
$$ \rho^T - (2m+2)\omega^T = \mathcal L_X \omega^T.$$
Hamiltonian holomorphic vector fields are defined
in Definition 4.4 and $\mathcal L_X$ stands for the Lie derivative by $X$. 
In the above equation the imaginary part of $X$ is necessarily a Killing
vector field. In general a Sasaki-Ricci soliton is a Sasaki-Einstein metric if and only if
$f_1$ vanishes. 

By definition a Sasaki manifold is 
toric if its K\"ahler cone is 
toric (see also Definition \ref{def42}). One of our  main results is stated as follows.

\begin{theorem}\label{MainThm} Let $S$ be a compact toric Sasaki manifold 
 with $c_1^B > 0$ and $c_1(D) = 0$.
Then there exists a Sasaki metric which is a 
Sasaki-Ricci soliton. In particular $S$ admits a Sasaki-Einstein metric
if and only if $f_1$ vanishes.
\end{theorem}

This is a Sasakian version of a result of Wang and Zhu (\cite{Wang-Zhu}) in the K\"ahlerian case.
As an application of Theorem \ref{MainThm} we obtain the following.

\begin{theorem}\label{thm2} Let $S$ be a compact toric Sasaki manifold
 with $c_1^B > 0$ and $c_1(D) = 0$. Then by
deforming the Sasaki structure varying  the Reeb field we get a Sasaki-Einstein structure.
\end{theorem}

Very recently, new obstructions to Sasaki-Einstein metrics
were studied by Gauntlett, Martelli, Sparks and Yau ( \cite{GMSY}).
Theorem 1.2 matches their expectation that those new obstructions are
not obstructions in the toric case, see section 6 of \cite{GMSY}.

A Sasaki manifold is said to be quasi-regular if all the leaves of the Reeb foliation
are compact, and irregular otherwise. A quasi-regular Sasaki manifold is said to be
regular if the Reeb foliation is obtained by a free $S^1$-action.

\begin{corollary}\label{cor2} There exists an irregular toric Sasaki-Einstein metric on 
the circle bundle of a power of the anti-canonical line bundle of the two-point blow-up
of the complex projective plane.
\end{corollary}

Examples of irregular Sasaki-Einstein metrics have been only recently known by
Gauntlett, Martelli, Sparks and Waldrum (\cite{GMSW1}, \cite{GMSW2}, \cite{MaSp1},
\cite{MaSp2}).
 They include an irregular toric Sasaki-Einstein metric on the circle
bundle associated with the canonical bundle of the one-point
blow-up of the complex projective plane. As far as the authors know an irregular example for
the two-point blow-up case has not been known. The irregularity of the example
mentioned in Corollary \ref{cor2} follows from the computations given by 
Martelli, Sparks and Yau (\cite{MSY1}); there is an earlier computation by 
Bertolini, Bigazzi and Cotrone for the dual theory through AdS/CFT correspondence
\cite{BeBiCo04}. 
For a different derivation of irregular Sasaki-Einstein metric on $S^2\times S^3$ see 
Hashimoto, Sakaguchi and Yasui \cite{HSY04}, and 
for other new Sasaki-Einstein metrics see also
\cite{new1} and \cite{new2}.

It should be mentioned that, before those irregular examples were found, many quasi-regular examples were constructed by Boyer, Galicki, Koll\'ar
and their collaborators. These constructions are surveyed in the article \cite{BG04} 
(see also \cite{BGK}).

This paper is organized as follows. In section 2 and 3 we review Sasaki geometry, and 
give a rigorous treatment for the transverse K\"ahler geometry. In section 4 we introduce
Hamiltonian holomorphic vector fields, and define the integral invariants $f_k$ as 
Lie algebra characters of the Lie algebra of all Hamiltonian holomorphic vector fields.
In section 5 we introduce Sasaki-Ricci soliton and a generalized integral
invariant which is useful for the study of Sasaki-Ricci soliton. We also set up
Monge-Amp\`ere equation on Sasaki manifolds to prove the existence of 
Sasaki-Ricci soliton.
In section 6 we review known facts about toric Sasaki manifolds. In particular we
describe the space of all K\"ahler cone metrics of Sasaki structures whose the basic K\"ahler class
is equal to $1/(2m+2)$ times basic first Chern
class. This space plays an important role in the study of the volume functional.
In section 7 we show that we can
give a $C^0$ estimate and completes the proof of the
existence of Sasaki-Ricci solitons on compact toric manifolds of  positive basic first Chern class. 
The point is that we can use the Guillemin metric to get a nice initial metric so that the analysis
of Wang and Zhu can be used.
In section 8 we will see that
on compact toric Sasaki manifolds we can always find a Sasaki structure with $f_1 = 0$,
thus proving Theorem \ref{thm2}. Then we prove Corollary \ref{cor2}.

The authors are grateful to Charles Boyer for pointing out our careless statements of 
the results without the condition $c_1(D) = 0$ in the earlier versions of the paper.

\section{Sasaki manifolds}

In this section we describe a Sasaki manifold in terms of two complex structures.
One is given on the cone over the Sasaki manifold, and the other is given
on the transverse holomorphic structure for a one dimensional foliation, called the
Reeb foliation. Proofs of the results in this section are concisely 
summarized in the papers of Boyer and Galicki \cite{boga99} and Martelli, Sparks and
Yau \cite{MSY2}. 

Let $(S, g)$ be a Riemannian manifold, $\nabla$ the Levi-Civita connection of
the Riemannian metric $g$, and let $R(X,Y)$ and $Ric$ respectively denote
the Riemann curvature tensor and Ricci tensor of $\nabla$.

\begin{definition}\ \ $(S,g)$ is said to be a Sasaki manifold if the cone manifold
 $(C(S ),\overline{g})=
(\R_+\times S , \ dr^2+r^2g)$ is K\"ahler.
\end{definition}

$S$ is often identified with the submanifold $\{r = 1\} = \{1\} \times S \subset C(S)$.
Thus the dimension of $S$ is odd and denoted by $2m + 1$. 
Hence $\dim_{\bfC}C(S) = m+1$. Let $J$ denote the
complex structure of $C(S)$. Define a vector field $\xi$ on $S$
and a $1$-form $\eta$ on $S$ 
by
\begin{equation}
\xi = J\frac{\partial}{\partial r},\ \ \eta(Y) = g(\xi, Y)
\end{equation}
where $Y$ is any smooth vector field on $S$. 
One can see that
\begin{enumerate}
\item $\xi$ is a Killing vector field on $S$;
\item the integral curve of $\xi$ is a geodesic;
\item $\eta (\xi)=1$ and $d\eta (\xi, X)=0$ for any vector field $X$ on $S$. 
\end{enumerate}
The vector field $\xi$ is called the
{\it characteristic vector field} or the {\it Reeb field}. 
The $1$-dimensional foliation generated by $\xi$ is called the {\it Reeb foliation}.
The 1-form
$\eta$ defines a $2m$-dimensional vector sub-bundle $D$ of the tangent bundle $TS$,
where at each point $p\in S$ the fiber $D_p$ of $D$ is given by
\begin{equation}\nonumber
D_p=\Ker\,\eta_p.
\end{equation}
There is an orthogonal decomposition of the tangent bundle $TS$
\[TS=D \oplus L_\xi,\]
where $L_\xi$ is the trivial bundle generated by the Reeb field
$\xi$. 

We next define a section $\Phi$ of the endomorphism bundle $\mathrm{End}(TS)$ of the tangent bundle $TS$
by setting $\Phi|_D = J|_D$ and $\Phi|_{L_\xi} = 0$ where we identified $S$ with
the submanifold $\{r=1\} \subset C(S)$. $\Phi$ satisfies 
\[\Phi^2=-I+\eta\otimes \xi\]
and
\[g(\Phi X, \Phi Y)=g(X, Y)-\eta(X)\eta(Y).\]
One can see that $\Phi$ can also be defined by
\begin{equation}
\Phi(X) = \nabla_X \xi
\end{equation}
for any smooth vector field $X$ on $S$. 
From these descriptions it is clear that $g|_D$ is an Hermitian metric on $D$.
An important property of Sasaki manifolds is that
$\frac12(d\eta)|_D$ is the associated 2-form of the Hermitian metric $g|_D$:
\begin{equation}
(d\eta)(X,Y) = 2 g(\Phi(X),Y)
\end{equation}
for all smooth vector fields $X$ and $Y$.
Hence $d\eta$
defines a symplectic form on $D$, and $\eta$ is 
a contact form on $S$ in the sense that $\eta \wedge (\frac12 d\eta)^m$
is nowhere vanishing. In particular
$\eta \wedge (\frac12 d\eta)^m$ defines a volume element on $S$.

We have $\xi$ and $\eta$ also on the cone $C(S)$ by putting
$$ \xi = Jr\frac{\p}{\p r}, \qquad  \eta(Y) = \frac 1{r^2}\overline{g}(\xi, Y) $$
where $Y$ is any smooth vector field on $C(S)$. Of course $\eta$ on $C(S)$ is 
the pull-back of $\eta$ on $S$ by the projection $C(S) \to S$.
Since  $L_\xi J=0$ the complex vector field $\xi-iJ\xi=\xi+ir\frac {\p}{\partial r}$ is a holomorphic vector field on $C(S)$, which preserves the cone structure. There is a holomorphic 
$\C^*$-action on $C(S)$
generated by $\xi-iJ\xi$. The local orbits of this action then give the Reeb foliation a transversely 
holomorphic structure.
As we will see in the next section this transversely holomorphic
foliation on $S$ is a K\"ahler foliation in the sense that it has a bundle-like transverse K\"ahler metric.

It is this transverse K\"ahler structure that we wish to study in this paper. We will study in later sections the deformations of transverse K\"ahler structures. 
They correspond to the deformations of $\eta$ fixing $\xi$. A fixed choice
of $\xi$ is called the polarization of the Sasaki manifold in \cite{BGS}. As $\eta$ varies the contact
bundle $D$ varies, while we fix the transverse holomorphic structure
of the Reeb foliation and the holomorphic structure of the cone $C(S)$, see Proposition 4.2
below.

To conclude this section we summarize the well-known facts about Sasaki geometry. First of all 
$\Phi$ defined above on the Sasaki manifold $S$ satisfies
\begin{equation}
(\nabla_X\Phi)(Y) = g(\xi,Y)X - g(X,Y)\xi.
\end{equation}
for any pair of vector fields $X$ and $Y$ on $S$.  

Conversely, if $\xi$ is a Killing vector field of unit length on a Riemannian manifold $S$ and 
the $(1,1)$-tensor $\Phi$ defined by $\Phi(X) = \nabla_X\xi$ satisfies (4), then the
cone $(C(S), dr^2 + r^2g)$ is a K\"ahler manifold and thus
$S$ becomes a Sasaki manifold. Here the complex structure $J$ on $C(S)$ is
defined as follows : 
$$ J r\frac{\partial}{\partial r} = \xi, \qquad JY = \Phi(Y) - \eta(Y) r\frac{\partial}{\partial r} \qquad \mathrm{for}\ Y \in \Gamma(TS).$$
Moreover the following conditions are equivalent and can be used as a definition of Sasaki
manifolds.
\begin{itemize}
\item [(2.a)] There exists a Killing vector field $\xi$ of unit length on $S$
so that the tensor field $\Phi$ of type $(1,1)$, defined by
$\Phi(X) ~=~ \nabla_X \xi$,
satisfies (4).
\item[(2.b)]
There exists a Killing vector field $\xi$ of unit length on $S $
so that the Riemann curvature satisfies the condition 
$$R(X,\xi)Y ~=~ g(\xi,Y)X-g(X,Y)\xi,$$
for any pair of vector fields $X$ and $Y$ on $S.$ 
\item[(2.c)] There exists a Killing vector field $\xi$ 
of unit length on $S $
so that the sectional curvature of every section containing $\xi$ equals one.
\item[(2.d)] The 
metric cone  $(C(S ),{\bar g})=
(\R_+\times\SS , \ dr^2+r^2g)$ over $S $ is K\"ahler.
\end{itemize}
Returning to the geometry of $C(S)$, we recall the following facts.
The $1$-form $\eta$ is expressed on $C(S)$ as 
\begin{equation}
\eta = 2d^c\log r
\end{equation}
where $d^c=\frac i2 (\bar \p-\p).$
This easily follows from 
$$\eta(Y) = \frac 1{r^2} \overline{g}(\xi, Y) = \frac 1{r^2} \overline{g}(Jr\frac{\p}{\p r}, Y). $$
From (5), the K\"ahler form of the cone $(C(S), dr^2 + r^2g)$ is expressed as
\begin{equation}
\frac 12 d(r^2\eta) = \frac12 dd^c r^2.
\end{equation}

\section{Transverse holomorphic structures and transverse 
K\"ahler structures}
In the sequel, we always assume that $\SS$ is a Sasaki manifold
with $(\xi, \eta, g, \Phi)$.
To understand the Sasaki structure well, 
we need to exploit the {\it transverse K\"ahler} structure on $S$.
Let ${\cal F}_\xi$ be the {\it Reeb foliation}
generated by $\xi$. 
As we saw in the previous section $\xi-iJ\xi=\xi+ir\frac {\p}{\partial r}$ is a holomorphic vector
field on $C(\SS)$, and there is an action on $C(\SS)$ of the holomorphic flow
generated by $\xi-iJ\xi$. The local orbits of this action defines a
transversely holomorphic structure on the Reeb foliation ${\cal F}_\xi$ in the following sense.

Let $\{U_\a\}_{\a\in A}$ be an open covering of $\SS$ and
$\pi_\a:U_\a\to V_\a\subset 
\C^m$ submersions such that when $U_\a\cap U_\b\neq \emptyset$
$$\pi_\a\circ \pi_\b^{-1}:\pi_\b(U_\a\cap U_\b)\to \pi_\a(U_\a\cap U_\b)$$
is biholomorphic. On each $V_\a$ we can give a K\"ahler structure as follows.
Let $D=\Ker\ \eta \subset T\SS$. There is a canonical isomorphism 
\[d\pi_\a:\DD_p \to T_{\pi_\a(p)} V_\a,\]
for any $p\in U_\a$. Since $\xi$ generates isometries  the restriction of the Sasaki
metric $g$ to $\DD$ gives 
a well-defined Hermitian 
metric $g^T_{\a}$ on $V_\a$. This Hermitian structure is in fact K\"ahler,
which can be seen as follows.

Let $z^1,z^2,\cdots z^m$ be the local holomorphic coordinates on $V_\a$. We pull
back these to $U_\a$ and still write them as $z^1,z^2,\cdots z^m$.
Let $x$ be the coordinate along the  leaves with $\xi=\frac{\p}{\p x}$.
Then $x,z^1,z^2,\cdots z^m$ form local coordinates on $U_\a$.
$(\DD\otimes \bfC)^{1,0}$ is spanned by the vectors of the form
\[\frac{\p}{\p z^i}+a_i\xi, \quad i=1,2\cdots,m.\]
It is clear that
\[a_i=-\eta(\frac{\p}{\p z^i}).\]
Since $i(\xi) d\eta =0,$
\[d\eta (\frac{\p}{\p z^i}+a_i\xi,\overline{\frac{\p}{\p z^j}+a_j\xi})
=d\eta (\frac{\p}{\p z^i}, \frac{\p}{\overline\p z^j}).\]
Thus the fundamental 2-form $\omega_\a$ of the Hermitian metric $g^T_{\a}$  on $V_\a$ is 
the same
as the restriction of $\frac12 d\eta$ to the slice $\{x= \mathrm{constant}\}$ in $U_\a$. Since the restriction
of a closed 2-form to a submanifold is closed in general, then $\omega_{\a}$ is 
closed. By this construction
$$\pi_{\a}\circ \pi_{\b}^{-1}:\pi_{\b}(U_{\a} \cap U_{\b})\to \pi_{\a}(U_{\a}\cap U_\b)$$
gives
an isometry of K\"ahler manifolds. Therefore, the foliation thus defined is a
transversely K\"ahler foliation. The collection of K\"ahler metrics $\{g^T_{\a}\}$ on $\{V_\a\}$ is 
called a transverse
K\"ahler metric. 
Since they are isometric over the overlaps we suppress $\a$ and denote by $g^T$.
We also write $\n^T$, $R^T$, $Ric^T$ and $s^T$ for
its Levi-Civita connection, the curvature, 
the Ricci tensor and the scalar curvature; 
of course these are collections of those defined on $\{V_\a\}$. 
It should be emphasized that, though $g^T$ are defined only locally on each $V_{\a}$, the
pull-back to $U_{\a}$ of the K\"ahler forms $\omega_{\a}$ on $V_{\a}$ patch together 
and coincide with the global form  $\frac12 d\eta$ on $S$, and $\frac12 d\eta$ can even be lifted to 
the cone $C(S)$ by the pull-back. For this reason we often refer to $\frac12 d\eta$ as the 
K\"ahler form of the transverse K\"ahler metric $g^T$. 
Although it is possible to re-define  $g^T$ at this stage as a global tensor on $S$ by setting
$ g^T = \frac 12 d\eta(\cdot, \Phi \cdot)$,
we do not take this view point because this may lead to a confusion about the space where
the K\"ahlerian geometry is performed. However, 
notice also that the transverse scalar curvature
 $s^T$ also lifts to $S$ as a global function which together with the lifted K\"ahler form
 $\frac 12 d\eta$ on $S$ is often used later in order to study the global properties
 of the Sasaki structures, e.g. to define integral invariants.
 
One can easily check the following. For $\wt X, \wt Y, \wt Z, \wt W \in \Gamma(D)$ 
and $X, Y, Z, W \in \Gamma(TV_\a)$
with $d\pi_\a(\wt X) = X, d\pi_\a(\wt Y) = Y, d\pi_\a(\wt Z) = Z, d\pi_\a(\wt W) = W$ we have
\begin{eqnarray}
\n^T_XY&=&d\pi_\a(\n_{\widetilde X} \widetilde Y),\\
\n_{\wt X}{\wt Y}&=&\wt {\n^T_XY} -g(JX,Y)\xi,\\
R(\wt X,\wt Y,\wt Z,\wt W)&=&R^T(X,Y,Z,W)-g(J\wt Y,\wt W)g(J\wt X,\wt Z)\nonumber \\
&& + g(J\wt X,\wt W)g(J\wt Y,\wt Z),\\
Ric^T(X,Z)&=&Ric(\wt X,\wt Z)+2g(\wt X,\wt Z).
\end{eqnarray}
\begin{definition} 
A Sasaki manifold 
$(\SS, g)$ is $\eta$-Einstein if there are two constants $\l$ and $\nu$ such that
\[Ric=\l g+\nu \eta\otimes \eta.\]
\end{definition}

In this case $\l + \nu = 2m$ always holds since $Ric(\xi,\xi) = 2m$ (see (2.b) of the
previous section).

\begin{definition} 
A Sasaki-Einstein manifold is a Sasaki manifold $(S,g)$with $Ric = 2m\, g$.
\end{definition}

Notice that a Sasaki manifold satisfying the Einstein condition is necessarily Ricci positive.

\begin{definition}
A Sasaki manifold $\SS$ is said to be transversely K\"ahler-Einstein
if
\[Ric^T=\tau g^T,\]
for some real constant $\tau$.\end{definition}

It is well-known that if $\SS$ is
a transversely K\"ahler-Einstein
Sasaki manifold if and only if $(\SS, g)$ is $\eta$-Einstein (cf. \cite{BGM}). In fact, if $Ric^T=\tau g^T$ then
\begin{equation}\label{Ricci1}
Ric = (\tau - 2)g + (2m + 2 - \tau) \eta \otimes \eta.
\end{equation}
Conversely if $Ric = \l g + \nu \eta \otimes \eta$ then
\begin{equation}
Ric^T = (\lambda + 2)g^T.
\end{equation}

Thus if $\l > -2$ then $Ric^T$ is positive definite. In this case by the $\cal D$-homothetic
transformation $g' = \a g + \a (\a -1) \eta\otimes\eta$ with $\a = \frac{\l +2}{2m + 2}$
we get a Sasaki-Einstein metric $g'$ with $Ric(g') = 2m\, g'$ (Tanno \cite{tanno}).

Thus there exists a Sasaki-Einstein metric if and only if there exists a transverse K\"ahler-Eisntein
metric of positive transverse Ricci curvature. In other words if we can find an obstruction to the
existence of transverse K\"ahler-Einstein metric of positive transverse Ricci tensor then  it is an
obstruction to the existence of Sasaki-Einstein metric. The invariant $f$ given in Theorem 
\ref{Futaki-inv} is one of such obstructions. The invariant $f$ has been also defined in 
\cite{BGS}.

On the other hand, the Gauss equation relating the curvature of submanifolds to 
the second fundamental form shows that {\it a Sasaki metric is Einstein if and only if
the cone metric on $C(S)$ is Ricci-flat K\"ahler}. In particular the K\"ahler cone of a Sasaki-Einstein
manifold has trivial canonical bundle. In section 8 we will reformulate $f$ as an obstruction
for a K\"ahler cone $C(S)$ with flat canonical bundle to admit a Ricci-flat K\"ahler cone metric.
From (\ref{Ricci1}) we see that $\tau = 2m+2$ for a Sasaki-Einstein metric, and thus we have
$$Ric = 2m g$$
and
\begin{equation}\label{Ricci2}
Ric^T = (2m+2)g^T.
\end{equation}
This also follows from (10).

\section{Basic forms and deformations of Sasaki structures}
Let $\SS$ be a compact Sasaki manifold.

\begin{definition} A $p$-form $\a$ on $\SS$ is called basic if
\[i(\xi)\a=0,\quad L_\xi \a=0.\]
Let $\Lambda_B^p$ be the sheaf of germs of basic $p$-forms and $\Omega_B^p
=\Gamma(\SS, \L_B^p)$ the set of all global sections of $\Lambda_B^p$.
\end{definition}

Let $(x, z^1,\cdots,z^m)$ be the coordinates system $U_\a$ given above.
We call such a coordinate system a
{\it foliation chart}. If $U_\a\cap U_\b \neq \emptyset$ and $(y, w^1,\cdots, w^m)$
is the foliation chart on $U_\b$, then
\[\frac {\p z^i}{\p \barw^j}=0, \quad \frac {\p z^i}{\p y}=0.\]
These mean that the form of type $(p,q)$
\[\a=a_{i_1 \cdots i_p\bar j_1\cdots \bar j_q} dz^{i_1}\wedge \cdots
dz^{i_p}\wedge d \bar{z}^{j_1} \wedge\cdots \wedge d \bar{z}^{j_q}\]
is also of type $(p,q)$ with respect to $(y, w^1, \cdots, w^m)$. 
If $\a$ is basic, then $a_{i_1 \cdots i_p\bar j_1\cdots \bar j_q}$ does not depend on $x$.
We thus have  the well-defined operators
\[\begin{array}{rcl}
\p_B:\Lambda_B^{p,q} &\to &  \Lambda_B^{p+1,q},\\
\barpartial_B:\Lambda_B^{p,q}& \to & \Lambda_B^{p,q+1}.\end{array}\]
It is easy to see that $d\a$ is basic if $\a$ is basic. If we set $d_B=d|_{\Omega^p_B}$ then
we have $d_B=\p_B+\bar \p_B$. Let $d^c_B=\frac i2 (\barpartial_B-\p_B).$
It is clear that
\[d_Bd^c_B=i\p_B\barpartial_B, \quad\quad d_B^2=(d^c_B)^2=0.\]
Let $d^*_B:\Omega^{p+1}_B\to \Omega^{p}_B$ be the adjoint operator 
of $d_B:\Omega^p_B\to \Omega^{p+1}_B$. The basic Laplacian $\D_B$
is defined 
\[\D_B=d^*_Bd_B+d_Bd_B^*.\]
Thus we can consider the {\it basic} de Rham complex $(\Omega^{\ast}_B, d_B)$ and
the {\it basic} Dolbeault complex $(\Omega^{p,\ast}, \barpartial_B)$ whose
cohomology groups are called the {\it basic} cohomology groups. Similarly, we
can consider the {\it basic} harmonic forms. Results of El Kacimi-Alaoui  (\cite{El}) assert
that we have the expected isomorphisms between basic cohomology groups and 
the space of basic harmonic forms.

Suppose $(\xi,\eta,\Phi,g)$ defines a Sasaki structure on $\SS$.
We define a new Sasaki structure fixing $\xi$ and varying $\eta$ as follows.
Let $\varphi \in \Omega^0_B$ be a smooth basic function. 
Put 
\[\tilde \eta= \eta + 2d^c_B\varphi.\]
It is clear that
\[d\tilde \eta=d\eta+2d_B d_B^c\varphi=d\eta+2i\p_B\bar \p_B \varphi.\]
For small $\varphi$, $\tilde\eta$ is non-degenerate in the sense that
$\tilde\eta\wedge (d\tilde\eta)^m$ is nowhere vanishing. 
\begin{proposition} Given a small smooth basic function $\varphi$, 
there exists a Sasaki structure on $S$ with the same $\xi$, the same holomorphic 
structure on the cone $C(S)$ and the same transversely
holomorphic structure of the Reeb foliation $\cal F_\xi$ but with the new contact
form $\wt \eta = \eta + 2d^c_B\varphi$.
\end{proposition}
\begin{proof}
Put
\begin{equation}\label{ato1}
\wt r = r\exp(\varphi).
\end{equation}
Then $\frac 12 dd^c \wt r^2$ gives a new K\"ahler structure on $C(S)$ and thus
a new Sasaki structure on $S$. Obviously the holomorphic structure on $C(S)$ is
unchanged for the new K\"ahler structure. 
By (\ref{ato1}),
$$ \wt \eta = 2 d^c\log \wt r = \eta + 2d^c_B\varphi. $$
It follows from this that the Reeb field is also unchanged 
because $\varphi$ is basic.
From the expression $(\ref{ato1})$ one sees
$$ \wt r \frac \p {\p \wt r} = r \frac \p {\p r}. $$
This shows that the transverse holomorphic structures of the
Reeb foliation is also unchanged.
\end{proof}

Thus the deformation 
$$\eta \to \tilde \eta=\eta+2d_B^c\varphi$$ gives
a deformation of Sasaki structure with the same
transversely holomorphic foliation and the same holomorphic structure
on the cone $C(S)$ and deforms the transverse K\"ahler
form in the same basic (1,1) class. We call this class the {\it basic K\"ahler class}
of the Sasaki manifold $S$.
Note that, however, the contact bundle $D$ may change under such a  deformation.

Now we define a $2$-form $\rho^T$ called the {\it transverse Ricci form} as follows.
This is first defined as a collection
of $(1,1)$ forms $\rho^T_\a$on $V_\a$ given by
$$ \rho^T_\a = - i \p\barpartial \log \det (g^T_\a).$$
These are just the Ricci forms of the transverse K\"ahler metrics $g^T_\a$. 
One can see that the pull backs $\pi_\a^{\ast} \rho^T_\a$ by 
$\pi_\a : U_\a \to V_\a$ patch together to give a global basic 2-form on $S$,
which is our  $\rho^T$. As in the K\"ahler case $\rho^T$ is $d_B$ closed and
define a basic cohomology class of type $(1,1)$. The basic cohomology class
$[\rho^T]$ is independent of the choice of the transverse K\"ahler form
$\frac 12 d\eta$ in the fixed basic $(1,1)$ class. The basic cohomology class $[\frac 1{2\pi}
\rho^T]$ is called the {\it basic first Chern class} of $S$, and is denoted by $c_1^B(S)$.
We say that the basic first Chern class of $S$ is positive (resp. negative) 
if $c_1^B(S)$ (resp. $- c_1^B(S)$) is represented by a basic K\"ahler form; This condition is
expressed by $c_1^B > 0$ (resp. $c_1^B < 0$). When $c_1^B > 0$ or $c_1^B < 0$ then
we say that $c_1^B$ is definite. 
If there exists a transverse K\"ahler-Einstein metric then the basic first Chern class
has to be positive, zero or negative according to the sign of the constant $\tau$
with $Ric^T = \tau g^T$. In particular if $S$ has a Sasaki-Einstein metric then the
basic first Chern class is positive. There are further necessary condition for the existence
of positive or negative transverse K\"ahler-Einstein metric:

\begin{proposition}\label{c_1(D)} The basic first Chern class is represented by
$\tau d\eta$ for some constant $\tau$ if and only if 
 $c_1(D) = 0$ where $D = \mathrm{Ker}\,\eta$ is the
contact bundle.
\end{proposition}
\begin{proof} This proposition follows from the long exact sequence (c.f. \cite{tondeur}):
 $$ \longrightarrow H^0_B(S) \overset{\delta}{\longrightarrow} H^2_B(S)
 \overset{i_2}{\longrightarrow} H^2(S;\bfR) \longrightarrow$$
 where $\delta a = a [d\eta]$ and $i_2 [\rho] = [\rho]$. Note $i_2 c_1^B =
 c_1(D)$ in general. 
\end{proof}
Thus if $S$ admits a transverse K\"ahler-Einstein metric 
$\rho^T = \tau \omega^T$ then $c_1(D) = 0$.
As was pointed out by Boyer, Galicki and Matzeu
\cite{BGM}, in the negative and zero basic first class case with $c_1(D) = 0$ the results of
El Kacimi-Alaoui \cite{El} together with Yau's estimate \cite{yau78} imply that
the existence of transverse K\"ahler-Einstein metric. One of main purposes of
the present paper is to consider the remaining positive case, assuming $c_1(D) = 0$.

\begin{pro} The quantity
\[\int_{\SS} \rho^T\wedge (d\eta)^{m-1}\wedge \eta\]
is independent of the choice of Sasaki structures with the same basic K\"ahler
class.\end{pro}

\begin{proof}Let $\varphi \in \Omega_B^0$ and put $\eta_t=\eta+td^c_B\varphi$
for small $t$. It is enough to show that
\[\frac d{dt}\int_{\SS}\rho^T_t\wedge(d\eta_t)^{m-1}\wedge \eta_t=0.\]
It is not difficult to check that
\begin{eqnarray*}
\frac d{dt}\int_\SS\rho^T_t\wedge(d\eta_t)^{m-1}\wedge \eta_t&=&
\int_\SS d_B d^c_B(\D_B\varphi)\wedge (d\eta_t)^{m-1}\wedge \eta_t\\
&&+(m-1)\int_\SS\rho_t^T \wedge (d\eta_t)^{m-2}\wedge d_Bd_B^c\varphi\wedge\eta_t)\\
&&
+\int_\SS\rho_t^T \wedge (d\eta_t)^{m-1}\wedge d^c_B\varphi).
\end{eqnarray*}
Since $\rho_t^T,$ $d\eta_t$, $d^c_B\varphi$ are basic, we have
$\rho_t^T \wedge (d\eta_t)^{m-1}\wedge d^c_B\varphi=0$.
It is clear that $d_Bd_B^c \D_B\varphi=dd^c_B \D_B\varphi$, from which we have
\[\int_\SS d_B d^c_B(\D_B\varphi)\wedge (d\eta_t)^{m-1}\wedge \eta_t
=-\int_Sd_B^c\D_B\varphi\wedge (d\eta_t)^m.\]
Since $d_B^c\D_B\varphi\wedge (d\eta_t)^m$ is a basic $(2m+1)$-form, it is zero.
Similarly, we have
\[ \int_\SS\rho_t^T \wedge (d\eta_t)^{m-2}\wedge d_Bd_B^c\varphi\wedge\eta_t)
=\int \rho_t^T \wedge (d\eta_t)^{m-2}\wedge d_B^c\varphi\wedge d\eta_t)=0.\]
This proves the proposition.
\end{proof}
This proposition means that the average $\overline{s}$ of the transverse scalar curvature $s^T$ 
depends only on the basic K\"ahler class where
\begin{equation}
\overline{s} = \frac{\int_S m\rho^T \wedge (\frac12 d\eta)^{m-1} \wedge \eta} {\int_S (\frac12 d\eta)^m \wedge \eta}
= \frac{\int_S s^T (\frac12 d\eta)^{m} \wedge \eta} {\int_S (\frac12 d\eta)^m \wedge \eta}.
\end{equation}

Recall, for the next definition, that $\frac12 d\eta$ is the transverse K\"ahler form. 
\begin{definition} A complex vector field $X$ on a Sasaki manifold
is called a Hamiltonian holomorphic vector field if
\begin{itemize}
\item [(1)] $d\pi_\a (X)$ is a 
 holomorphic vector field on $V_\a$;
\item [(2)] the complex valued function $u_X:=\sqrt{-1}\eta(X)$ satisfies
\[\bar \p_B u_X=-\frac{\sqrt{-1}}{2} i(X) d\eta.\]
\end{itemize}
Such a function $u_X$ is called a Hamiltonian function.
\end{definition}

If $(x, z^1,\cdots, z^m)$ is a foliation chart on $U_\a$, then $X$ is written
as
\[X= \eta(X)\frac{\p}{\p x}+\sum_{i=1}^m X^i\frac{\p}{\p z^i}-
\eta(\sum_{i=1}^m X^i\frac{\p}{\p z^i})\frac{\p}{\p x},\]
where $X^i$ are local holomorphic basic functions.
 
 Note that $X+ i \eta(X) r\frac{\p}{\p r}$ is a holomorphic 
 vector field on ${\cal C}(\SS)$.
 A Hamiltonian holomorphic vector field $ X$ is the orthogonal projection of
 a Hamiltonian holomorphic vector field $\wt X$ on $C(\SS)$ to $S = \{r=1\}$, whose 
 Hamiltonian function $\tilde u$ satisfies $\xi(\tilde u)=\frac{\p}{\p r} \tilde
 u=0$, ie., $\tilde u$ is basic and homogenous of degree zero with respect to $r$.
 
 \begin{remark}
 If $u_X = c (= \mathrm{constant})$ then by $\mathrm{(2)}$ of Definition 4.4, $X = c\xi$ and 
 $d\pi_\a (X) = 0.$ In the case of regular Sasakian manifolds this corresponds to the 
 fact that the Hamiltonian function 
 for a Hamiltonian vector field on a  symplectic manifold is unique up to constant.
\end{remark}
 
 Let $\hh$ denote the set of all  Hamiltonian holomorphic 
 vector fields. One can easily check that $\hh$ is a Lie algebra.
 Nishikawa and Tondeur \cite{NT} proved that if the scalar curvature $s^T$ of the
 transverse K\"ahler metric is constant then $\mathfrak h$ is reductive, extending
 Lichnerowicz-Matsushima theorem in the K\"ahler case.

Let $\omega_\a$ be the K\"ahler metric on $V_\a$, {\it i.e.}
\[\pi^*_\a \omega_\a =\frac 12 d\eta.\]
Then as is well-known, there is a smooth function $f_\a$ on $V_\a$
 such that
\[\omega_\a=i\p \bar \p f_\a.\]
Hence in terms of the foliation chart
$(x, z^1,\cdots, z^m)$ on $U_\a$,
\[\eta =dx + 2 d^c_B f_\a\]
and
\begin{equation}
d\eta=2d_Bd^c_B f_\a=2i\p_B\bar \p_B f_\a.\end{equation}

\begin{pro} By the deformation $\eta \to \tilde \eta=
\eta+2d^c_B\phi$, $u_X$ is deformed to
$u_{X}+X\phi$.\end{pro}
\begin{proof} Note that $u_X= i \eta(X)$ and $\tilde u_{X}= i\tilde \eta(X).$
\end{proof}
\begin{pro}\label{LichOp}A complex valued basic function $u$ is a
Hamiltonian function for some Hamiltonian holomorphic vector field $X$ if and only if
\begin{equation}
\D_B^2u+(i\p_B\bar\p_Bu, \rho^T)+(\bar \p_B u,\bar \p_B s^T)=0\end{equation}
where $\rho^T$ and $s^T$ are the transverse Ricci form and transverse scalar curvature.
\end{pro}
\begin{proof} From $\bar \p_B u=- \frac i2\,i(X)d\eta$, we have
$\frac{\p u}{\p \barz^j}=g^T_{i\bar j} X^i$,
ie. $X^i=(g^T)^{i\bar j}\frac{\p u}{\p \barz^j}$. 
For simplicity of the notation, we omit $T$ and $B$ in the proof.
$X^i=(g^T)^{i\bar j}\frac{\p u}{\p \barz^j}$ is holomorphic if and only if
\[\begin{array}{ll}
 \Leftrightarrow & \ds\vs \n_{\bari} \n_{\bar j}u=0, \quad\forall i,\,j\\
 \Leftrightarrow & \ds\vs \n^{\bari}\n^{\barj}\n_{\bari} \n_{\bar j}u=0,\\
\Leftrightarrow &\ds\vs  \D^2u+R^{i\bar j}u_{i\bar j}+ u^is_i=0\\
\Leftrightarrow & \ds \D^2u+(i\p\bar \p u,\rho)+ (\bar\p u, \bar \p  s)=0
\end{array}\]
\end{proof}

The following result was obtained also by Boyer, Galicki and Simanca \cite{BGS}.
\begin{theorem}\label{Futaki-inv}
 Let $\eta_t=\eta+td^c_B\phi$ be the deformation
 of Sasaki structures. Let $X$ be a Hamiltonian holomorphic vector
 field on $\SS$ and $u_t$ the Hamiltonian functions with respects to $X$ and $\eta_t$.
 Then
 \[f(X)=-\int_\SS u_t(m\rho^T_t\wedge (\frac12 d\eta_t)^{m-1}\wedge \eta_t- \bar s (\frac12 d\eta_t)^m\wedge \eta_t)
 \]
 is independent of $t$, where $\overline{s}$ was difined in $\mathrm{(15)}$.
 \end{theorem}

\begin{proof} A direct computation using Proposition 4.6 gives
\begin{eqnarray}
- \frac d{dt} f(X,t) &=& \int_\SS(\bar \p_B u, \bar\p_B \varphi)
(m\rho^T_t\wedge (\frac12 d\eta_t)^{m-1}\wedge \eta_t- \bar s (\frac12 d\eta_t)^m\wedge \eta_t)\label{a1}\\
&&+\int_\SS umd_Bd_B^c(\D_B\varphi)\wedge (\frac12 d\eta_t)^{m-1}\wedge \eta_t \label{a2}\\
&&+\int_\SS um\rho^T_t\wedge (m-1)(\frac12 d\eta_t)^{m-2}\wedge d_B d_B^c\varphi\wedge \eta_t \label{a3}\\
&&+\int_\SS u(m\rho^T_t\wedge (\frac12 d\eta_t)^{m-1}-\bar s(\frac12 d\eta_t)^m)\wedge d_B^c\varphi \label{a4}\\
&&+\int u\bar s \D_B\varphi (\frac12 d\eta_t)^m\wedge \eta_t\label{a5}\end{eqnarray}
Here, and later too, we suppressed the suffix $t$ in the notations of inner product and the basic Laplacian.
It is easy to check that (\ref{a5}) + the second term of (\ref{a1}) vanishes.
(\ref{a4}) vanishes, for its integrand is basic of degree $(2m+1)$. (\ref{a3}) equals 
to
\[\int_\SS \varphi d_Bd_B^c u\wedge m\rho^T_t\wedge (m-1)(\frac12 d\eta_t)^{m-2}\wedge\eta_t.\]
Thus,
\begin{eqnarray}
\label{f'}
\ds \vs - \frac d{dt} f(X,t) &=& \ds
\int_\SS (\bar \p_B u, \bar\p_B \varphi)\,s^T(\frac12 d\eta_t)^m\wedge \eta_t-\int_\SS u\,\D_B^2\varphi\,
(\frac12 d\eta_t)^m\wedge\eta  \\
&&+\ds \int \varphi\, d_Bd_B^c u\wedge m\rho^T\wedge (m-1)(\frac12 d\eta_t)^{m-2}\wedge\eta_t.\nonumber
\end{eqnarray}
Take any point $p \in U_\a \subset S$ and a foliation chart $(x, z^1, \cdots, z^m)$ on $U_\a$ such that,
on $V_\a$, $\p/\p z^1,\ \cdots, \p/\p z^m$ are orthonormal and that either $(\partial_i\p_\barj u)$ or
$(R^T_{i\barj})$ is diagonal. Then the second term on the right hand side of (\ref{f'}) is equal to
\begin{equation*}
\int_S \varphi\,\sum_{I\ne j} u_{i\bari} R^T_{j\barj}\,(\frac12 d\eta)^m \wedge \eta =
\int_S \varphi\, (-\D_Bu\,s^T - (d_Bd_B^c  u, \rho^T))(\frac12 d\eta)^m \wedge \eta.
\end{equation*}
From this and Proposition \ref{LichOp} we finally get
\begin{equation*}
\ds \vs - \frac d{dt} f(X,t) = \ds
- \int_\SS \varphi(\D_B^2u + (\bar \p_B u, \bar\p_B s^T) + (d_Bd_B^c u, \rho^T))(\frac12 d\eta)^m\wedge\eta
= 0.
\end{equation*}
\end{proof}

The linear function $f$ on the Lie algebra $\mathfrak h$ of all Hamiltonian holomorphic vector fields 
is obviously an obstruction to the existence of transverse
K\"ahler metric  of constant scalar curvature
in the fixed basic K\"ahler class. In particular it obstructs the existence of
transverse K\"ahler-Einstein metric, extending earlier result of the first author for 
Fano manifolds,  see \cite{F}. 

The invariant $f$ has different expressions.
In the positive  case, we assume that $c_1^B= (2m+2)[\frac12 d\eta]_B$
which is assured by the assumption $c_1(D) = 0$, see Proposition \ref{c_1(D)}.
By a result of El Kacimi-Alaoui
\cite{El}, there is a basic function $h$ such that
\begin{equation}\label{EL}
\rho^T- (2m+2)\frac12 d\eta =i \p_B\bar \p_B h.\end{equation}
In this case, the average of the scalar curvature is $m$. Thus we have,
by the definition of $f$, that
\begin{equation}\label{f}
\begin{array}{rcl}
f(X)&=&\ds\vs -\int_\SS u_X m(\rho\wedge (\frac12 d\eta)^{m-1}\wedge \eta-
(2m+2) (\frac12 d\eta)^{m}\wedge \eta)\\
&=&\ds\vs -m\int _\SS u_X i\p_B\bar\p_B h\wedge (\frac12 d\eta)^{m-1}\wedge \eta\\
&=& \ds - \int_S u_X \D_Bh (\frac12 d\eta)^{m}\wedge \eta
= \int_S  Xh\ (\frac12 d\eta)^{m}\wedge \eta.
\end{array}\end{equation}
 From  (\ref{f}) it is clear that when the Sasaki manifold $\SS$ has a transverse K\"ahler-Einstein metric 
 in its basic K\"ahler class, then $f$ vanishes.
 
We also have the following generalization.
Recall that $\Omega$ is the Sasaki cone, namely
\[\Omega=\{\tilde\eta= \eta+ 2d_B^c\varphi\,|\, \varphi \hbox{ is basic and } 
(d\tilde\eta)^m\wedge \tilde\eta \hbox{ is nowhere vanishing } \}.\]
Let $\o_B=\frac12 d\eta$ be the transverse K\"ahler form and $\Theta$ be the curvature matrix of 
the transverse Levi-Civita connection of $\omega_B$. Consider
the basic $2m$-form
\[c_m(\o_B\otimes I+\frac{\sqrt{-1}}{2\pi} \lambda\Theta),\]
where $I$ denotes the identity matrix, $c_m$ the invariant polynomial corresponding
to the determinant. We can expand it as follows
\[c_m(\o_B\otimes I+\frac{\sqrt{-1}}{2\pi} \lambda\Theta)
=\o_B^m+\l c_1(w_B)\wedge \o_B^{m-1}+\cdots+\l^mc_m(\o_B),\]
where $c_i(\o_B)$ is the $i$-th Chern form with respect to $\o_B$. 
It is clear that a tangential vector field to $\Omega$ at $\eta$ 
can be expressed by a basic function $\psi$ with
a normalization
\[\int_\SS \psi (\frac12 d\eta)^m\wedge\eta=0.\]
Hence the space of all basic functions is the tangent space to $\Omega$.
We define a one-form $\a$ on $\Omega$ by 
\begin{equation}
\label{aa}
\psi\to\int_\SS\phi c_m( \o_B\otimes I+\frac{\sqrt{-1}}{2\pi} \lambda\Theta)
\wedge\eta.
\end{equation}

Define $f_k:\hh\to {\mathbb C}$
by
\[f_k(X)=\int_\SS u_X c_k(\o_B)\wedge (\frac12 d\eta)^{m-k}\wedge \eta,\]
where $c_k(\o_B)$ is the $k$th Chern form defined above and $u_X$ is the
Hamiltonian function for $X$ with the normalization
\[\int_\SS u_X (\frac12 d\eta)^m\wedge\eta=0.\]

\begin{theorem}\label{thm28} $f_k$ is independent of the choice of $\eta$.
Moreover $f_k$ is a Lie algebra homomorphism and
obstructs the existence of a Sasaki metric of harmonic basic $k$th Chern form.
Moreover $f_1 = 1/2\pi\,f$.
\end{theorem}

\begin{proof}The proof is similar to the proof of Theorem \ref{Futaki-inv}. 
See Appendix.
\end{proof}
\begin{lem}\label{lemma26}
The one-form $\a$ defined by $\mathrm{(\ref{aa})}$ is a closed form on $\O$.\end{lem}
\begin{proof} To show the Lemma, it is convenient to consider an arbitrary map
 $\gamma:[0,1]\times[0,1] \to\Omega$. Let $(s,t)$ be the coordinates on 
 $[0,1]\times[0,1]$.
 By definition, we know that 
 \[(\gamma^*\a)(\frac \p {\p s})=
 \a (\gamma_* \frac{\p}{ \p s})=\a(\frac {\p \gamma} {\p s}).\]
Thus,
\[\begin{array}{rcl} \gamma^*\a &=&\ds\vs
\left(\int_\SS \frac {\p \gamma} {\p s} c_m(\o_B\otimes I+\frac {\sqrt {-1}}{2\pi} \l 
\Theta)\wedge \eta\right)ds\\
&&\ds+\left(\int_\SS \frac {\p \gamma} {\p t} c_m(\o_B\otimes I+\frac {\sqrt {-1}}{2\pi} \l 
\Theta)\wedge \eta\right)dt,\end{array}\]
and $\a$ is closed if only if $\gamma^*\a$ is closed for any $\gamma$. 
In view of Lemma \ref{lemma24} and Lemma \ref{lemma25},
the latter can be proved as in the K\"ahler (see \cite{F}).\end{proof}

\begin{theorem}\label{thm30}
Let $\eta$ and $\eta'$ be two Sasaki structures in $\Omega$ and 
$\eta_t=\eta+ 2d^c_B\varphi_t$ $(t\in [a,b])$ be a path in $\Omega$ connecting $\eta$ and $\eta'$. Then
\[{\cal M}_k(\eta,\eta')=
\int_a^b\int_\SS \dot \varphi_t(c_k(\o_B)-H(c_k(\o_B)))\wedge (\frac12 d\eta)^{m-k}
\wedge \eta\]
is independent of the path $\eta_t$, where
$H(c_k(\o_B))$ is the basic harmonic part of $c_k(\o_B)$ with respect to
the transverse K\"ahler form $\omega_B = \frac12 d\eta$.
\end{theorem}

\begin{proposition}\label{pro32}
Fix $\eta_0\in \Omega$ and consider the functional
$m_k:\Omega \to {\mathbb R}$ defined by
\[m_k(\eta)=-{\cal M}_k(\eta_0,\eta).\]
then the critical points of $m_k$ are the 
 Sasaki metrics of harmonic basic $k$-th Chern form.
 When $k=1$ and $[c_1(\frac12 d\eta)]_B=const.[\frac12 d\eta]_B$, then the critical points 
 are the metrics of transverse K\"ahler-Einstein. 
\end{proposition}

The proofs of Theorem \ref{thm30} and Proposition \ref{pro32} can be given by the principle stated in
the Appendix. 
Theorem \ref{thm28} and Theorem \ref{thm30} respectively extend results of Bando \cite{bando83} and 
Bando and Mabuchi \cite{bandomabuchi86}.

\section{Sasaki-Ricci solitons}
To study the existence of Sasaki-Einstein metrics, or equivalently
transversely K\"ahler-Einstein metrics of positive scalar curvature (or $\eta$-Einstein metric),
a natural analogue of K\"ahler Ricci flow
was introduced in \cite{SW}. 
Assume that $c_1^B=\kappa [\frac12 d\eta]_B$, where $\kappa$ is normalized so that
$\kappa=-1,0,$ or $1$ for simplicity.
We consider the following flow $(\xi, \eta(t),\Phi(t),g(t))$ with 
initial data $(\xi, \eta(0),\Phi(0),g(0))=(\xi,\eta,\Phi,g)$
\begin{equation}\label{flow1}
\frac{d}{dt}g^T(t)=-(Ric^T_{g(t)}-\kappa g^T(t)),\end{equation}
or equivalently
\begin{equation}\label{flow2}
\frac{d}{dt}  d\eta(t)=-(2\rho^T_{g(t)}-\kappa d\eta(t)).\end{equation}
This flow is called {\it Sasaki-Ricci flow}.
 {Locally}, if we write $\frac12 d\eta=\sqrt{-1}g^T_{i\bar j}dz^i\wedge
d \bar z^j$, we can  check that
$\rho^T=-\sqrt{-1}\partial \bar \partial \log \det
(g^T_{k\bar l })$. Let $\eta(t)=\eta+2d_B^c\varphi(t)$, a family of basic functions
$\varphi(t)$
Then the  flow can be written as
\begin{equation}\label{flow3}
\frac d{dt} \varphi =\log 
 \det (g_{i\bar j}^T+\varphi_{i\bar j})
 -\log (\det g_{i\bar j}^T)+\kappa \varphi-h,\end{equation}
 where $h$ is a basic function defined by
 \begin{equation}\label {el} 
 \rho^T_g- \frac{\kappa}2 d\eta = d_Bd_B^c h. 
 \end{equation}
 The solvability of (\ref{el}) was proved in \cite{El}.
 The well-posedness of the flow was
proved in \cite{SW}. Like the K\"ahler-Ricci flow \cite{cao},
the long-time existence can be also proved. When the flow converges, then the limit
is a transverse K\"ahler-Einstein metric. In fact, one can show that
when $\kappa=-1,$ or $0$, then the flow globally converges to an
$\eta$-Einstein metric. See also \cite{El} and \cite{BGM}.
Hence, the remaining interesting case is when $\kappa=1$, namely the basic first Chern 
form of the Sasaki manifold is of positive definite. But from now on, we assume
that $\kappa=2m+2$ because this normalization fits to the study of Sasaki-Einstein
metric, see (\ref{Ricci2}).
In this case,
 in general the convergence of the flow could not be obtained. What one can hope
 is the limit converges in some sense to a soliton solution, as in the K\"ahler
 case. A Sasaki structure $(\SS, \xi,\eta,\Phi,g)$ with a Hamiltonian
 holomorphic vector field $X$ is called a {\it transverse K\"ahler-Ricci soliton} or 
 {\it Sasaki-Ricci soliton} if
 \begin{equation}\label{soliton1}
 Ric^T-(2m+2)g^T=\mathcal L_Xg^T,\end{equation}
 or equivalently,
 \begin{equation}\label{soliton2}
 2\rho^T-(2m+2)d\eta=\mathcal L_X(d\eta).\end{equation}
 
 In the next section, we will prove that on any toric Sasaki manifold
 there always exists a Sasaki-Ricci soliton. To end this section, we give a
 generalization of the invariant $f$, which is an obstruction of
 the existence of the Sasaki-Ricci solitons.
 
  Recall that there is a basic function satisfying
 (\ref{EL}). As in \cite{F}, we define the following operator 
 \[\begin{array}{rcl}
 \ds\vs \D_{B}^h u = \D_B u-\n^i u \n_i h,\end{array}
 \] where $\n=\n^T$ is the Levi-Civita connection of the transverse K\"ahler metric.
 One can show that the operator  $\D_{B}^h$ is self-adjoint in the following
 sense
 \begin{equation}\label{*}
 \begin{array}{rcl}
\ds \vs  \int_\SS \D_{B}^h u\ \overline{v}\ e^h (\frac12 d\eta)^m\wedge \eta 
 &=& \ds \int_\SS u\  \overline{\D_{B}^h v}\ e^h (\frac12 d\eta)^m\wedge \eta \\
 &=&\ds \int_\SS (\bar \p_B u, \bar \p_B v)\ e^h (\frac12 d\eta)^m\wedge \eta.\end{array}
 \end{equation}

 We need to consider  ``{\it normalized} Hamiltonian 
 holomorphic vector fields", whose corresponding
 Hamiltonian functions $u_X$ satisfying
 \begin{equation}\label{normal}
  \int_\SS u_X e^h (\frac12 d\eta)^m\wedge\eta =0.\end{equation}
For  any Hamiltonian 
 holomorphic vector field $X$, there is a unique constant $c\in {\mathbb R}$
 such that $X+c\xi$ is a normalized Hamiltonian 
 holomorphic vector field. For simplicity of notation, from now on 
any  holomorphic vector field $X$ we consider is normalized and its Hamiltonian function
is denoted by $\t_X$. Hence $\t_X$ satisfies (\ref{normal}).

 The operator  $\D_{B}^h$ has  the following properties, whose proof
 can be given as in \cite{F}.
 
 \begin{theorem}\label{thm47} We have
 \begin{itemize}
 \item [(1)] The first eigenvalue $\l_1$ of $\D_B^h$ is greater than or equal to $2m+2$.
 \item[(2)] The equality $\l_1=2m+2$ holds if and only if there exists  a non-zero Hamiltonian 
 holomorphic vector field $X$.
 \item [(3)] $\{u \in \Omega_B(\SS)^{\mathbb C}\,|\, \D_B^h u=(2m+2)u\}$ is isomorphic to
 $\{X\,|\,$normalized Hamiltonian 
 holomorphic vector fields$\}$. The correspondence is given, in a local foliation chart,
 by
 \[u\to u\xi+\n^i u\frac \p{\p z^i}- \eta(\n^i u\frac \p{\p z^i})\xi.\]
 \end{itemize}\end{theorem}
 
 \begin{remark}
 \begin{enumerate}
 \item[(1)] $\n^i \t=g_T^{i\bar j} \frac {\p u}{\p z^i}$ is holomorphic.\\
 \item[(2)] Let $u\in \{u \in \Omega_B(\SS)^{\mathbb C}\,|\, \D_B^h u=(2m+2)u\}$.
  By putting $v=1$ in (\ref{*}), we get
 \[\int_\SS u e^h (\frac12 d\eta)^m\wedge \eta=0\]
 \end{enumerate}
 \end{remark}

 Now, as in Tian and Zhu \cite{TG02} we define a generalized invariant
$f_X$  
 for a given Hamiltonian holomorphic vector field $X$ by
 \[ f_X(v) =-\int_\SS  \theta_v e^{\theta_X}(\frac12 d\eta)^m\wedge \eta.\]
 
 We will leave the proof of the invariance of $f_X$ to the reader. It is also
 easy to check that $f_X$ is an obstruction of the existence of Sasaki-Ricci
 solitons as follows.
 
 Let $(\SS, g, X)$ be a Sasaki-Ricci soliton, i.e., we have (\ref{soliton2}).
  From the above discussion we know that
  \[dd_B^c h=L_X(\frac12 d\eta)=d(i_X ( \frac12 d\eta))= id\bar \p_B \t_X.\]
  It follows that $\t_X=h + \mathrm{constant}$. Hence $f_X(v)=0$ for all $v\in {\mathfrak h}$.
  The next proposition shows that we can always find such an $X \in {\mathfrak h}(\SS)$.
  
\begin{proposition}\label{pro49} There exists an $X\in {\mathfrak h}(\SS)$ such that
\[f_X(v)=0, \quad \forall v\in {\mathfrak h}(\SS).\]
\end{proposition}

\begin{proof} The proof can be given by arguments similar to \cite{TG02},
and we will not reproduce them here.
\end{proof}

Now we wish to set up the Monge-Amp\`ere equation to prove the existence of a transverse
K\"ahler-Ricci soliton for the choice of $X$ in Proposition \ref{pro49}. Choose an initial
Sasaki metric $g$ such that the transverse K\"ahler form $\omega^T = \frac 12 d\eta$ represents
the basic first Chern class of the normal bundle of the Reeb foliation. There exists a
smooth basic function $h$ such that
\begin{equation}\label{SR1}
\rho^T - (2m+2)\omega^T = i\p_B\barpartial_B h.
\end{equation}
Suppose we can get a new Sasaki metric $\tildeg$
satisfying the Sasaki-Ricci soliton equation by a transverse K\"ahler deformation.
Let 
$$\tildeomega^T = i\, \tildeg^T_{i\barj}\, dz^i \wedge d\barz^j
=  i\, (g^T_{i\barj} + \varphi_{i\barj})\, dz^i \wedge d\barz^j, $$
$\widetilde{\rho}^T$ and $\widetilde{\theta}_X$ respectively denote the transverse K\"ahler 
form, transverse Ricci form and the Hamiltonian function for the
normalized Hamiltonian function $X$ with respect to $\tildeg$
where $\varphi$ is a smooth basic function $S$. Then by the Sasaki-Ricci  
soliton equation we have
\begin{equation}\label{SR2}
\widetilde{\rho}^T - (2m+2)\tildeomega^T = L_X \tildeomega^T = i \p_B \barpartial_B \widetilde{\theta}_X.
\end{equation}
As one can see easily (c.f. Appendix 2, \cite{futakimabuchi95})
\begin{equation}\label{SR3}
\widetilde{\theta}_X = \theta_X + X\varphi.
 \end{equation}
From (\ref{SR1}), (\ref{SR2}) and (\ref{SR3}) we get
\begin{equation}\label{SR4}
\frac{\det (g^T_{i\barj} + \varphi_{i\barj})}{\det (g^T_{i\barj})} = \exp( - (2m+2)\varphi - 
\theta_X - X\varphi + h)
\end{equation}
with $ (g^T_{i\barj} + \varphi_{i\barj}) $ positive definite (recall that $\varphi$ is a basic 
function).
In order to prove the existence of a solution to (\ref{SR4}) we consider a family of equations
parametrized by $t \in [0,1]$:
\begin{equation}\label{SR5}
\frac{\det (g^T_{i\barj} + \varphi_{i\barj})}{\det (g^T_{i\barj})} = \exp( - t (2m+2)\varphi - \theta_X - X\varphi + h)
\end{equation}
with $ (g^T_{i\barj} + \varphi_{i\barj})$ positive definite. It is sufficient to show that the subset of
$[0,1]$ consisting of all $t$ for which (\ref{SR5}) has a solution is non-empty, open and closed.
Combining the arguments of \cite{El}, \cite{yau78} and \cite{Zhu} one can show that (\ref{SR5}) 
has a solution at $t=0$ and that the openness is also satisfied by the implicit function theorem. By 
El Kacimi-Alaoui's generalization of Yau's
estimates (\cite{yau78}) for the transverse Monge-Amp\`re equations it suffices to show the $C^0$ estimate for $\varphi$ 
to prove to the closedness.
To get the $C^0$ estimate for $\varphi$ it is sufficient to get the $C^0$ estimate of on an
open dense subset of $S$. We prove in section 7 that this $C^0$ estimate can be obtained
for toric Sasaki manifolds.

\section{Toric Sasaki manifolds}
In this section, we recall known facts about toric Sasaki manifolds,
following \cite{Lerman}, \cite{BG0}, and a slight modification of some arguments of \cite{MSY1}.

\begin{definition}\label{def34} A co-oriented contact structure $D$ on a manifold
$M$ is a codimension $1$ distribution such that
\begin{itemize} 
\item[(1)] $D^0=\{\a\in T^*M\,|\, \a(X)=0,\, \forall X\in D\}$ is an oriented real line
bundle, and a component $D^{0}_+$ of $D^0\slash\{0\}$ is chosen.
\item[(2)] $D^0\slash\{0\}$ is a symplectic submanifold of $T^*M$.
\end{itemize}\end{definition}
Note that $T^*M$ has a natural symplectic structure $\sum_{i=1}^n dq^i\wedge dp_i$
where $q^1, \cdots, q^n$ is local coordinates on $M$ and 
$p = \sum_{i=1}^np_idq^i$ is a cotangent vector.

\begin{definition}\label{def36} Suppose a Lie group $G$ acts on $M$ preserving
$D$, and consequently $D^0$. The moment map for the action of $G$ is a map
$\mu:D^0_+\to \g^*$ defined by
\[\la \mu(q,p),X\ra=\la p, X_M(q)\ra,\]
where $X_M(q):= \left.\frac d{dt}\right|_{t=0}exp (tX)\cdot q$ is a vector field on $M$ induced by 
$X\in\g$. Here $\g$ is Lie algebra of $G$ and $\g^*$ is its dual.
 \end{definition}
\begin{definition}\label{def37}
\ \ A contact toric $G$-manifold is a co-oriented contact manifold $(M,D)$ with an action of 
a torus $G$ preserving $D$ and with $2\dim G =dim M+1$.\end{definition}
 From \cite{Lerman}, we have
 \begin{lemma}\label{lem38} Let $(M,D)$ be a toric contact manifold with an action of a 
 torus $G$. Then zero is not in the image of the contact moment map $\mu: D^0_+\to \g^*$.
 Moreover the  moment cone $C(\mu)$
 defined  by
 \[C(\mu):= \mu(D^0_+)\cup \{0\}\]
 is a convex rational polyhedral cone.
 \end{lemma}
 
 Recall that a subset ${\cal C}\subset \g^*$ is a convex rational polyhedral set, if
 there exists a finite set of vectors $\{\l_j\}$ in the integral lattice
 ${\mathbb Z}_G:=\Ker \{\exp :\g\to G\}$ and $\mu_i\in {\mathbb R}$ such that 
 \[{\cal C} =\bigcap_j\{X\in \g^*\,|\, \la \l_j, X \ra + \mu_i \ge 0\}.\]
Clearly it is a cone if all $\mu_i$ are 0.

 Let $\eta$ be a contact form, i.e., $\Ker\, \eta =D$ and $d\eta|_{D}$ non-degenerate. Let $G$ 
 be a torus action on $M$ preserving $\eta$. The moment map $\mu_\eta:M\to \g^*$ is defined by
 \[\la \mu_\eta(x), X\ra=\eta( X_M(x)), \quad \forall x\in M.\]
 The contact form $\eta$ is a section $\eta: M\to D^0_+$ and we obviously have
 \[\eta^*\mu=\mu_{\eta},\]
 since $\la \eta^*\mu, X\ra=\eta (X_M).$ For a $G$-invariant contact form $\eta$
 we define the moment cone $C(\mu_\eta) $ by
 \[C(\mu_\eta)=\{r e\,|\, e\in \mu_\eta(M), r\in [0,\infty)\}.\]
 It is clear that
 \[C(\mu_\eta)=C(\mu).\]
 
 Now we introduce the notion of toric Sasaki manifolds.
 \begin{definition}\label{def42} A toric Sasaki manifold
 $\SS$ is a Sasaki manifold of dimension $2m+1$ with Sasaki structure
 $(\xi,\eta,\Phi, g)$ such that there is an effective action 
 of $(m+1)$-dimensional torus $G$ preserving the  Sasaki structure
 and that $\xi$ is an element of the Lie algebra $\g$ of $G$. Equivalently, a toric Sasaki manifold
 is a Sasaki manifold whose K\"ahler cone is a toric K\"ahler manifold.
 \end{definition}
 
 \begin{proposition}\label{pro43}
 Let $\SS$ be a toric Sasaki manifold with Sasaki structure
 $(\xi,\eta,\Phi, g)$ and $\eta_t=\eta+2d_B^c\varphi$ a $G$-invariant basic deformation of 
 Sasaki structure where $\varphi$ is a $G$-invariant basic smooth function. Then the moment cones 
 $C(\mu_{\eta_t})$ are the same for all $t$, i.e.,
 \[C(\mu_{\eta_t})=C(\mu_{\eta}), \quad \forall t.\]
 \end{proposition}
 \begin{proof}This follows since every moment cone is rational polyhedral cone.
 \end{proof}
  
Let $G^c\cong ({\mathbb C}^*)^{m+1}$ denote the complexification of $G$. Then $G^c$
acts on the cone $C(\SS)$ as biholomorphic automorphisms. The moment map on $C(\SS)$
with respect to the K\"ahler form $\o=d(\frac12 r^2 \eta)$ is given by
\[\begin{array}{rcl}
\ds\vs \mu: C(\SS) &\to & \ds \g^*\\
\ds \la \mu(x), X\ra &=&\ds r^2\eta (X_\SS(x)),\end{array}\]
where we have used the natural diffeomorphism $C(\SS)\cong {\mathbb R}^+\times \SS$ and
view  vector fields on $\SS$ as vector fields on $C(\SS)$. 
Notice that we deleted $1/2$ so that there is a consistency with the moment maps for
contact manifolds.
It is clear that the image of $\mu$
is the same with the moment cone defined above, which is denoted by
$C(\mu)$. Let $\mathrm{Int}C(\mu)$ denote the interior of
$C(\mu)$. Then the action of $G$ on $\mu^{-1}(\mathrm{Int}C(\mu))$  is free and the orbit space
is  $\mathrm{Int}C(\mu)$. This means that $\mu^{-1}( \mathrm{Int}C(\mu))$
is a torus bundle over $\mathrm{Int}C(\mu)$. On the other hand the image 
$\mathrm{Im}(\mu_\eta)$ 
of the moment map $\mu _\eta:\SS\to \g^*$ is given by
\[\mathrm{Im}(\mu_\eta)=\{\a \in C(\mu)\,|\, \a(\xi)=1 \}.\]
The hyperplane $\{\a \in  \g^*\,|\, \a(\xi)=1\}$ is called
{\it characteristic plane} in \cite{BG0}. Notice that the constants differs
by $1/2$ from \cite{MSY1} because we use the moment map as a contact manifold.

In the rest of this section we study the Guillemin metric obtained by the
K\"ahler reduction through the Delzant construction (c.f. \cite{Abreu}, \cite{Guill}, \cite{Lerman}).

Assume that the moment cone $C(\mu)$ of our Sasaki manifold $S$ is described by
\begin{equation}\label{mcone}
C(\mu) = \{y \in \frak g^{\ast}\ | \ l_i(y) = \lambda_i \cdot y \ge 0,\ i=1, \cdots, d\},
\end{equation}
and let $C(\mu)^{\ast}$ be its dual cone
\begin{equation}\label{dualcone}
 C(\mu)^{\ast} = \{\tildex \in \frak g\ | \ \tildex\cdot y \ge 0\ \mathrm{for\ all}\ y \in C(\mu)\}. 
 \end{equation}
Then the Reeb field $\xi$ is considered as an element of the interior of $C(\mu)^{\ast}$
since $\frac12 r^2\eta(\xi) = \frac12 r^2 > 0$. 
We identify $\frak g^{\ast} \cong \bfR^{m+1} \cong \frak g$, and regard
$$ \lambda_i = (\lambda_i^1, \cdots, \lambda_i^{m+1}), \quad \xi = (\xi^1, \cdots, \xi^{m+1}).$$
As was shown in \cite{MSY1} the symplectic potential $G^{can}$ of
the canonical metric in the above sense is expressed by
$$ G^{can} = \frac 12 \sum_{i=1}^d l_i(y)\log l_i(y). $$
If we put 
$$ G_{\xi} = \frac 12 l_{\xi}(y) \log l_{\xi}(y) - \frac 12 l_{\infty}(y) \log l_{\infty}(y)$$
where
$$ l_{\xi}(y) = \xi\cdot y, \quad l_{\infty}(y) = \sum_{i=1}^d \lambda_i \cdot y,$$
then 
\begin{eqnarray}\label{sympotential1}
G_{\xi}^{can} &=& G^{can} + G_{\xi} \\
&=& \frac 12 \sum_{i=1}^d l_i(y) \log l_i(y) 
+ \frac 12 l_{\xi}(y) \log l_{\xi}(y) - \frac 12 l_{\infty}(y) \log l_{\infty}(y) \nonumber
\end{eqnarray}
gives a symplectic potential of a K\"ahler metric on $C(S)$ such that the induced
Sasaki structure on $S$ has $\xi$ as the Reeb field. To see this, compute
\begin{equation}\label{G1}
(G_{\xi}^{can})_{ij} = \frac12 \sum_{k=1}^d \frac{\lambda_k^i\lambda_k^j }{l_k(y)}
+ \frac12 \frac{\xi^i\xi^j }{l_{\xi}(y)} - \frac12 
\frac{\sum_{k=1}^d \lambda_k^i\sum_{\ell=1}^d\lambda_{\ell}^j }{l_{\infty}(y)}
\end{equation}
and observe that $\xi^i = 2\sum_j (G_{\xi}^{can})_{ij}y_j$ and that 
$\sum_{i,j}(G_{\xi}^{can})_{ij} y_iy_j = \frac12 l_{\xi}(y) > 0$, hence 
$(G_{\xi}^{can})_{ij}$ is positive definite.

Since any two complex structures associated to
a polytope are equivariantly biholomorphic (c.f. Proposition A.1 in \cite{Abreu})
we may assume that the complex structure obtained by the Delzant construction is the same
with the complex structure of the K\"ahler cone $C(S)$ of the Sasaki manifold $S$ under
consideration. Thus the Sasaki structure induced by the above Delzant construction has
the same complex structure and Reeb field. If we denote by 
$$\overline{\wt g} = d\wt r^2 + \wt r^2 \wt g$$ 
the Sasaki metric  of this Sasaki structure then we have
$$ \wt r \frac {\p}{\p \wt r} = J\xi = r \frac{\p}{\p r}.$$
This implies that $\wt r = r \exp(\varphi)$ for some basic smooth function $\varphi$.
Taking $d^c\log$ we get 
$$ \wt \eta = \eta + 2d^c \varphi.$$
This is a transverse K\"ahler deformation described in Proposition 4.2.

Thus we have proved the following.

\begin{proposition}\label{initial} Let $S$ be a compact toric Sasaki manifold and $C(S)$ its
K\"ahler cone. Let $\xi$ be the Reeb field.
Then we may assume that there is a transverse K\"ahler deformation of the
Sasaki structure of $S$ whose symplectic potential is of the form $\mathrm{(\ref{sympotential1})}$.
\end{proposition}

Now we assume hereafter that the initial Sasaki structure is so chosen that the 
symplectic potential $G$ is written as (\ref{sympotential1}). Let 
$$\wt x^j = \frac{\p G}{\p y_j} = \frac{\p G^{can}_{\xi}}{\p y_j}$$
be the inverse Legendre transform of $y_j$. Then 
$$(z^1, \cdots, z^{m+1}) = (\wt x^1 + i\wt \phi^1, \cdots,
\wt x^{m+1} + i \wt \phi^{m+1})$$ 
is the affine logarithmic coordinate system on
for $\mu^{-1}(\mathrm{Int}C(\mu)) \cong (\bfC^{\ast})^{m+1}$, i.e. the standard holomorphic 
coordinates are given by
$$(e^{\wt x^1 + i\wt \phi^1}, \cdots,
e^{\wt x^{m+1} + i \wt \phi^{m+1}}).$$
(We reserve the notation $x^j + i \phi^j$ for a subtorus $H^c$ of $G^c \cong (\bfC^{\ast})^{m+1}$
which will appear later.)

Suppose that the basic first Chern class $c_1^B$ of the normal bundle of the Reeb foliation is positive and that $c_1(D) = 0$. Then
the basic K\"ahler class $[\omega^T] = \frac 12 [d\eta]$ can be so chosen  that $c_1^B = (2m+2)[\omega^T]$.
If we denote by $F$ the K\"ahler potential on $C(S)$, $F$ and $G=G^{can}_{\xi}$ are related by
$$ F = \wt x\cdot y - G.$$
Then the Ricci form $\rho$ on $C(S)$ must be written as
\begin{equation}\label{F1}
 \rho = - i\p \barpartial \log \det (F_{ij}) = i \p \barpartial h
 \end{equation}
with 
$$ r\frac{\p}{\p r}h = 0, \qquad \xi h = 0,$$
i.e. $h$ is a pull-back of a basic smooth function $S$. This of course is equivalent to
$$ \rho^T = (2m+2)\frac 12 d\eta + i\p_B\barpartial_B h.$$
Since $T^{m+1}$-invariant pluriharmonic harmonic function is an affine function we see from
(\ref{F1}) that there exist $\gamma_1, \cdots, \gamma_{m+1} \in \bfR$ such that,
replacing $h$ by $h + \mathrm{const}$ if necessary,
\begin{equation}\label{F2}
\log \det (F_{ij}) = - 2 \gamma_i \wt x^i - h.
\end{equation}
In terms of $G$, (\ref{F2}) can be written as
\begin{equation}\label{F3}
\det(G_{ij}) = \exp(2\sum_{i=1}^{m+1}\gamma_i G_i + h).
\end{equation}

One can compute the right hand side of (\ref{F3})
using (\ref{sympotential1}) to get
\begin{equation}\label{F5}
\det(G_{ij}) = \Pi_j\left(\frac{l_j(y)}{l_{\infty}(y)}\right)^{(\lambda_j,\gamma)}
(l_{\xi}(y))^{-(m+1)}\exp(h).
\end{equation}
On the other hand using (\ref{G1}) we can compute the right hand side of
(\ref{F3}) to get
\begin{equation}\label{F6}
\det(G_{ij}) = f(y) \Pi_j (l_j(y))^{-1},
\end{equation}
where $f$ is a smooth positive function on $C(\mu)$.
It follows from (\ref{F5}) and (\ref{F6}) that
\begin{equation}\label{F7}
(\lambda_j,\gamma) = -1,\qquad j= 1,\cdots,d.
\end{equation}
Since the cone $C(\mu)^{\ast}$ is a cone over a finite polytope
there are $(m+1)$ linearly independent $\lambda_i$'s. Thus $\gamma$ is
uniquely determined from the moment cone $C(\mu)$, and is rational.
The equalities (\ref{F7}) show that the structure of the cone is very special. 
If $\gamma$ is a primitive lattice vector then the apex is
a Gorenstein singularity as Martelli, Sparks and Yau pointed out in section 2.2
of \cite{MSY2}.

Recall from \cite{MSY1} that
\begin{equation}\label{Euler}
r\frac{\p}{\p r} = 2 \sum_{j=1}^{m+1} y_j \frac{\p}{\p y_j},
\end{equation}
\begin{equation}\label{Reeb}
\xi_i = 2\sum_jG_{ij}y_j.
\end{equation}
The left hand side of (\ref{F3}) is homogeneous of degree $-(m+1)$ by (\ref{G1})
while applying
$ \sum_{j=1}^{m+1} y_j \frac{\p}{\p y_j}$ to the right hand side of (\ref{F3})
gives
$$ \sum_{j=1}^{m+1} y_j \frac{\p}{\p y_j}\exp(2\sum_{i=1}^{m+1}
\gamma_iG_i +h)
= (\xi, \gamma)\exp(2\sum_{i=1}^{m+1}\gamma_iG_i + h)$$
where we used (\ref{Reeb}). Hence we obtain
\begin{equation}\label{F4}
(\xi, \gamma) = -(m+1).
\end{equation}

The condition (\ref{F1}) says that the Hermitian metric $e^h\det(F_{ij})$ gives a flat
metric on the canonical bundle $K_{C(S)}$. Consider a holomorphic $(m+1)$-form
$\Omega$ of the form
$$ \Omega = e^{i\alpha}e^{\frac h2} \det(F_{ij})^{\frac12} dz^1 \wedge \cdots \wedge dz^{m+1}.
$$
From (\ref{F2}) $\Omega$ is written as
$$ \Omega = e^{i\alpha} \exp(-\sum_{i=1}^{m+1}\gamma_i\tildex^i)dz^1 \wedge \cdots \wedge
dz^{m+1}.$$
Hence taking $- \alpha = \sum_{i=1}^{m+1}\gamma_i\wt\phi^i$ we have
\begin{equation}\label{F8}
 \Omega = e^{-\sum_{i=1}^{m+1}\gamma_iz^i} dz^1 \wedge \cdots \wedge dz^{m+1}.
 \end{equation}
 Note that when $\gamma$ is not integral but only rational then $\Omega$ is multi-valued.
 Let $\ell$ be a positive integer such that $\ell\gamma$ is a primitive element of the integer
 lattice. Then $\Omega^{\otimes\ell}$ defines a section of $K^{\otimes\ell}_{C(S)}|_U$ where $U$ is the open
 dense subset of $C(S)$ corresponding to the interior of the moment cone. But since 
 $||\Omega^{\otimes\ell}|| = 1$, $\Omega^{\otimes\ell}$ extends to a smooth section of
 $K^{\otimes\ell}_{C(S)}$.
 
Since $\xi$ is expressed as $\xi = \sum_i\xi^i \frac{\p}{\p \wt\phi^i}$ we have
\begin{equation}\label{F9}
 \mathcal L_{\xi} \Omega = - i(\xi,\gamma)\Omega = i(m+1)\Omega.
 \end{equation}
From (\ref{F8}) and (\ref{F2}) we see
\begin{equation}\label{e112}
\left(\frac{i}{2}\right)^{m+1}(-1)^{m(m+1)/2}\Omega\wedge \overline{\Omega}=
\exp (h)\frac{1}{(m+1)!}\omega^{m+1}
\end{equation}
where $\omega$ denotes the K\"ahler form of $C(S)$.

To sum up we have obtained the following.
\begin{proposition}\label{tangent}
There exist a unique rational vector $\gamma \in \frak g^{\ast}$ such that (\ref{F7}) holds
and a multi-valued holomorphic $(m+1,0)$-form $\Omega$ with the following properties.
For some positive integer $\ell$, $\Omega^{\otimes\ell}$ defines a holomorphic section
of $K^{\otimes\ell}_{C(S)}$. 
Further for any K\"ahler cone
metric on $C(S)$ such that the Ricci form $\rho$ is written as
\begin{equation}\label{rho}
 \rho = i\p\barpartial h
 \end{equation}
where $h$ is the pull-back of a smooth basic function on $S$,
\begin{enumerate}
\item[(1)]
the Reeb field $\xi$ satisfies (\ref{F4}) and (\ref{F9});
\item[(2)]
if we denote by $\omega$ the K\"ahler form of the K\"ahler cone metric then
the equation (\ref{e112}) is satisfied.
\end{enumerate}

Conversely if the Reeb field $\xi$ satisfies either (\ref{F4}) or (\ref{F9}) then 
the K\"ahler cone metric satisfies (\ref{rho}) for some basic function $h$ on $S$.
\end{proposition}

\section{Sasaki-Ricci soliton on toric Sasaki manifolds}

In this section we want to show the existence of Sasaki-Ricci solitons on
any toric Sasaki manifolds with $c_1^B > 0$ 
as stated in Theorem \ref{MainThm} in the introduction.
As was explained in section 5 we have only to give a $C^0$ estimate for the family
of Monge-Amp\`ere equations (\ref{SR5}):

\begin{lemma}\label{C^0} Let $S$ be a compact toric Sasaki manifold
with $c_1^B > 0$ and $c_1(D) = 0$. Then
there exists a constant $C > 0$ independent of $t \in [0,1]$ 
such that $\sup_S |\varphi| \le C$ for any solution for $\mathrm{(\ref{SR5})}$.
\end{lemma}

The rest of this section is devoted to the proof of Lemma \ref{C^0}, which will completes
the proof of Theorem \ref{MainThm}.

First let us take any subtorus $H\subset G$ of codimension 1 
such that its Lie algebra ${\mathfrak h}$
does not contain $\xi$. Let $H^c \cong ({\mathbb C}^*)^m$ denote
the complexification
of $H$. Take any point $p\in \mu^{-1}(\mathrm{Int}C(\mu))$ and consider the orbit
$Orb_{C(S)}(H^c,p)$ of the $H^c$-action on $C(S)$ through $p$. 
Since $H^c$-action preserves $-J\xi=r\p\slash \p r$, it descends to an action on the set
$\{r=1\}\subset C(\SS)$, which we identify with the Sasaki manifold $\SS$. 
More precisely this action is described as follows. Let $\gamma : H^c \times C(S) \to
C(S)$ denote the $H^c$-action on $C(S)$.  Let ${\overline p}$ and $\overline{\gamma(g,p)}$ respectively
be the points on $\{r=1\}$ at which the flow lines through $p$ and $\gamma(g,p)$ generated by  $r\p\slash
\p r$ respectively meet $\{r=1\}$. 
Then the $H^c$-action on $S \cong \{r=1\}$
is given by $\bar{\gamma} : H^c \times \{r=1\} \to \{r=1\}$ where
$$\bar{\gamma}(g, {\overline p}) = \overline{\gamma(g,p)}.$$
Let $Orb_S(H^c,{\overline p})$ be the orbit of the induced action of $H^c$
on $\{r=1\} \cong \SS$.

\begin{proposition}\label{pro44}
The transverse K\"ahler structure of the Sasaki manifold $S$ over the open dense subset 
$\mu_{\eta}^{-1}(\mathrm{Int}\mathrm{Im}(\mu_{\eta}))$, which is the inverse image by
$\mu_{\eta}$ of the interior of the intersection of the characteristic hyperplane and 
$C(\mu)$, 
is completely determined by the restriction of $\frac 12 d\eta$ to $Orb_{C(S)}(H^c, p)$.
\end{proposition}
\begin{proof} Let $q \in Orb_{C(S)}(H^c, p)\subset C(\SS)$ be any point.
Since $r\p \slash 
\p r-i\xi=r\p \slash 
\p r-J(r\p \slash 
\p r)$ is preserved by $H^c$, a neighborhood $V_q$ of $q$ in $Orb_{C(S)}(H^c, p)$ is mapped biholomorphically to some $V_\a$ where $\barq\in U_\a \subset S$ and
$\pi_\a:U_\a\to V_\a$ is given by the transversely holomorphic structure of the Reeb foliation. Thus the transverse K\"ahler
structure on $V_\a$ is determined by $d\eta|_{V_q}$ because $(V_q, \frac12 d\eta|_{V_q})$ is isometric to
$(V_\a, \frac12 d\eta|_{V_\a})$ as K\"ahler manifolds; this is because 
$\eta = 2d^c \log r$ is homogeneous
of degree $0$.
But for any $
r\in \mu_\eta^{-1}(\mathrm{Int}\mathrm{Im}(\mu_{\eta}))$ the trajectory through $r$ generated
 by $\xi$ meets $Orb_S(H^c, \barp)$ and
$\xi$ generates a one parameter subgroup of isometries. So, the transverse K\"ahler geometry at
any $r$ is determined by the transverse K\"ahler geometry along the points on $Orb_S(H^c, \barp)$. 
This trajectory may meet $Orb_S(H^c, \barp)$ infinitely many times when the Sasaki structure is
irregular. But the transverse structures at all of them define the same
K\"ahler structure because $\xi$ generates a subtorus in $T^{m+1}$ and we assumed that $T^{m+1}$ preserves the Sasaki structure. 
Thus 
$\frac 12 d\eta|_{Orb_{C(S)}(H^c,p)}$ completely determines the transverse
 K\"ahler structure of $\SS$ on $\mu_\eta^{-1}
(\mathrm{Int}\mathrm{Im}(\mu_{\eta}))$.
\end{proof}

$(Orb_{C(S)}(H^c, p), \frac12 d\eta|_{Orb_{C(S)}(H^c, p)})$ and 
$(Orb_S(H^c, \barp), \frac12 d\eta|_{Orb_S(H^c, \barp)})$ are essentially the same 
in that if we give them the holomorphic structures induced from the holomorphic structure
of $H^c$ then they are isometric K\"ahler manifolds. The difference between them is that
$Orb_{C(S)}(H^c, p)$  is a complex submanifold of the complex submanifold $C(S)$
while $Orb_S(H^c, \barp)$ is a complex submanifold in the real Sasaki manifold $S$.
In what follows we are interested in the K\"ahler potential of $\frac 12 d\eta$. For this purpose
$Orb_{C(S)}(H^c, p)$ is better to treat.

On $Orb_{C(S)}(H^c, p) \cong  ({\mathbb C}^*)^m$ we use the affine logarithm coordinates 
$$(w^1,w^2,\cdots ,w^m)=
(x^1+\sqrt {-1}\t^1,x^2+\sqrt {-1}\t^2,\cdots, x^m+\sqrt {-1}\t^m)$$
for a point
$$(e^{x^1+\sqrt {-1}\t^1}, e^{x^2+\sqrt {-1}\t^2}, \cdots ,e^{x^m+\sqrt {-1}\t^m}) \in 
 ({\mathbb C}^*)^m \cong H^c.$$
Since $\eta$ is $H$-invariant, so is $d\eta$. Therefore $\frac12 d\eta|_{Orb_{C(S)}(H^c, p)}$ is determined by a convex 
function $u^0$ on ${\mathbb R}^m$, namely
\[\begin{array}{rcl}
\ds\vs \frac 12 d\eta|_{Orb_{C(S)}(H^c, p)} &=&\ds \sqrt{-1} \p\bar \p u^0\\
&=&\ds \frac{\sqrt{-1}}4
\frac{\p^2u^0}{\p x^i\p x^j} dw^i\wedge d \overline{w^j}.\end{array}\]
Hence we have
\begin{equation}\label{m1}
\begin{array}{rcl}
\ds
i(\frac \p{\p \t^i})\frac 12 d\eta|_{Orb_{C(S)}(H^c, p)} &=& \ds -\frac 12 \frac{\p^2u^0}{\p x^i\p x^j} dx^j
=\ds  -\frac 12 \,d(\frac {\p u^0}{\p x^i}).\end{array}\end{equation}
0n the other hand since $L_{\frac {\p}{\p \t^i}}\eta =0$, we have 
\[ i(\frac \p{\p \t^i})\frac12 d\eta =-\frac12 d (\eta(\frac \p{\p \t^i}))\]
and from this and (\ref{m1}) it follows that
\begin{equation}\label{m2}
\eta(\frac \p{\p \t^i})= \frac {\p u^0}{\p x^i}+c_i,
\end{equation}
where $c_i\in{\mathbb R}$ is a constant.

Now we wish to know more about the K\"ahler potential $u^0$. One way of expressing $u^0$ is
\begin{equation}\label{potential1}
u^0 = \log r|_{Orb_{C(S)}(H^c,p)} + \mathrm{const}.
\end{equation}
For, since $Orb_{C(S)}(H^c,p)$ is a complex submanifold of $C(S)$ and $\eta = 2d^c\log r$,
$$ (\frac 12 d\eta)|_{Orb_{C(S)}(H^c,p)} = (dd^c \log r )|_{Orb_{C(S)}(H^c,p)}
= dd^c(\log r |_{Orb_{C(S)}(H^c,p)}).$$
If we take the K\"ahler metric on $C(S)$ as the canonical metric obtained by the
K\"ahler reduction through the Delzant construction, 
we get a Sasaki structure for which the transverse K\"ahler potential $u^0$ has a more explicit description as explained below.

By (\ref{sympotential1}) the real part $\tildex^i$ of the affine
logarithm coordinates in $C(S)$ is given by
\begin{eqnarray}\label{x^j}
&&\tildex^j = \frac{\p G^{can}_{\xi}}{\p y_j}\\
&=& \frac 12 \sum_{i=1}^d \lambda^j_i (1+\log l_i(y)) + \frac 12 \xi^j(1 + \log l_{\xi}(y)) - 
\frac 12 \sum_{i=1}^d \lambda_i^j(1+\log l_{\infty}(y) )\nonumber \\
&=& \frac 12 \sum_{i=1}^d \lambda_i^j \log l_i(y) + \frac 12 \xi^j(1 + \log l_{\xi}(y)) -
\frac 12 \sum_{i=1}^d \lambda_i^j \log l_{\infty}(y).\nonumber
\end{eqnarray}
The K\"ahler potential $F_{\xi}^{can}$ on $C(S)$ is then obtained by the Legendre transform:
\begin{equation}\label{sympotential2}
F_{\xi}^{can}(\tildex) = \tildex\cdot y - G_{\xi}^{can}(y) = \frac 12 l_{\xi}(y).
\end{equation}
Now we know that $\frac{r^2}{2}$ is also a K\"ahler poteintial on $C(S)$, and hence
$\frac{r^2}{2} - \frac 12 l_{\xi}(y)$ is a harmonic function on $\bfR^{m+1}$. 
But $\frac 12 l_{\xi}(y)$ is 
bouded from below as (\ref{esti11}) is satisfied with $\bfu = \xi$, and $\frac{r^2}{2}$ is
also clearly bounded from below. Thus  $\frac{r^2}{2} - \frac 12 l_{\xi}(y)$ must be a constant.
But $\frac{r^2}{2} - \frac 12 l_{\xi}(y)$ tends to $0$ as $r$ tends to $0$. Therefore
$$ F_{\xi}^{can} = \frac 12 l_{\xi}(y) = \frac{r^2}{2}.$$
It follows that
\begin{eqnarray}\label{sympotential3}
u^0 &=& \log r|_{Orb_{C(S)}(H^c,p)} + \hbox{const}  \\
&=& \frac 12 \log (2F_{\xi}^{can})|_{Orb_{C(S)}(H^c,p)} + \hbox{const} \nonumber \\
&=& \frac 12 \log l_{\xi}|_{Orb_{C(S)}(H^c,p)} + \hbox{const}. \nonumber
\end{eqnarray}

Now we consider the moment map on the K\"ahler manifold 
$Orb_{C(S)}(H^c, p) \cong Orb_{S}(H^c, \barp)$ for the action of $H \cong T^m$. 
This is defined as 
\[j^*\circ \mu_\eta: Orb_{C(S)}(H^c, p) \to {\mathfrak h}^*\]
where $j:{\mathfrak h}\to \g$ is the inclusion and
$$\la \mu_\eta(y),X\ra = \eta(X)(y), \qquad X \in \mathfrak h,\ y \in Orb_{C(S)}(H^c, p),$$
$X$ being identified with a vector field on $Orb_{C(S)}(H^c, p)$.  This is essentially the same
as the restriction of the moment map $\mu_\eta : S \to \mathfrak h^*$ to $Orb_{S}(H^c, \barp)$.
Hence the image of $j^*\circ \mu_\eta$ is
equal to
\[j^*(\mathrm{Im}(\mu_\eta))=\{j^*\a\,|\, \a\in C(\mu), \a(\xi)=1 \}.\]
This is a (possibly irrational) compact convex polyhedron. Identifying $\mathfrak h$ with
$\bfR^m$ in the canonical way, the interior $\mathrm{Int}j^*\mathrm{Im}(\mu_{\eta})$ 
of $j^*\mathrm{Im}(\mu_{\eta})$ coincides up to translation with
\begin{equation}\label{Sigma}
\Sigma := \{ Du^0(x) = 
(\frac{\p u^0}{\p x^1}(x),\,\cdots,\,\frac{\p u^0}{\p x^m}(x))\,|\,x \in \bfR^m\}
\end{equation}
because of (\ref{m2}). Let $p_1,\,\cdots,\,p_\ell$ be the vertices of the closure 
$\overline{\Sigma}$ of $\Sigma$.

\begin{proposition}\label{pro46} Consider the Sasaki structure defined by the K\"ahler
cone metric on $C(S)$ with the symplectic potential (\ref{sympotential1}).
Let $u^0$ be the K\"ahler potential of $(Orb_{C(S)}(H^c,p), \frac 12 d\eta)$. 
Define $\overline{v} : \bfR^m \to \bfR$ by
$$ \barv(x) = \max_{1\le i \le \ell} \la p_i,x\ra. $$
Then there exists a constant $C$ such that $|u^0 - \barv| \le C$.
\end{proposition}
\begin{proof} From (\ref{sympotential3}) it is sufficient to show
$$ |2\barv(x) - \log l_{\xi}(y)| \le C$$
where $x$ and $y$ are related as follows: first write $\tildex$ for $jx \in \bfR^{m+1}$ 
take the Legendre transform
$$ y_j = \frac{\p F^{can}_{\xi}}{\p \tildex^j}$$
and then restrict $y_j$ to $Orb_{C(S)}(H^c,p)$.
For any vertex $q_i$ of $\{\a\in C(\mu)\,|\,\a(\xi) = 2\}$ such that
$$ j^{\ast}q_i = 2p_i, \qquad i = 1,\,\cdots,\,\ell$$
we have
$$ \la 2p_i,x\ra = \la j^{\ast}q_i,x\ra = \la q_i, jx \ra = q_i\cdot \tildex$$
where we again wrote $\tildex$ for $jx$ in the last term (we do so throughout the rest of the proof
of Proposition \ref{pro46}). Then from (\ref{x^j}) we have
\begin{equation}\label{esti1}
q_j\cdot \tildex - \log l_{\xi}(y) = \frac 12 \sum_{i = 1}^d l_i(q_j) \log l_i(y) -
\frac 12 l_{\infty}(q_j) \log l_{\infty}(y) + 1.
\end{equation}
In this proof we use the following simple fact repeatedly: {\it Let $\bfu$ be a non-zero
vector in $\bfR^{m+1}$ and $V$ be a closed strictly convex polyhedral cone in the 
open half space $\{y \in \bfR^{m+1}\,|
\,\bfu\cdot y > 0\} $. Then there are positive constants $c$ and $C$ such that for any 
$y \in V$ we have 
\begin{equation}\label{esti11}
 c|y| \le \bfu\cdot y \le C|y|. 
 \end{equation}
 }
These constants $c$ and $C$ will appear many times and take different values, but
we will use the same notation by taking smaller value of $c$ and lager value of $C$. 
This will not cause any problem as they appear only finitely many times.
Recall that $\sum_{i=1}^d \lambda_i$ is the Reeb field for the canonical metric (\cite{MSY1}).
Since $\sum_{i=1}^d \lambda_i$ is in the interior of $C(\mu)^{\ast}$ we have from the
above fact that for any $y \in C(\mu)$ 
\begin{equation}\label{l_infty}
c|y| \le l_{\infty}(y) \le C|y|.
\end{equation}
On the other hand by the Schwarz inequality we have for each $i$
$$ l_i(y) \le C|y|. $$
Let $q_1,\,\cdots,\,q_{\ell}$ be the vertices of $\{y \in C(\mu)\,|\, \ell_{\xi}(y) = 2\}
= \{\a\in C(\mu)\,|\,\a(\xi) = 2\}$. Suppose 
$$q_j \in \bigcap_{k=1}^m L_{i_k}$$
where $L_i$ denotes the hyperplane $\{l_i(y) = 0\}$. It follows from (\ref{esti1}) that
\begin{eqnarray}\label{esti2}
q_j\cdot \tildex - \log l_{\xi}(y) &=& \frac 12 \sum_{i \notin \{i_1,\cdots,i_m\}} l_i(q_j) \log l_i(y) \\
&&\quad- \frac 12 \sum_{i \notin \{i_1,\cdots,i_m\}} l_i(q_j) \log l_{\infty}(y) + 1\nonumber \\
&\le& \frac 12 \sum_{i \notin \{i_1,\cdots,i_m\}} l_i(q_j) (\log |y| + \log C)\nonumber\\
&&\quad- \frac 12 \sum_{i \notin \{i_1,\cdots,i_m\}} l_i(q_j) (\log |y| + \log c) + 1\nonumber \\
&\le& C.\nonumber
\end{eqnarray}
This proves
\begin{equation}\label{above}
 2\barv(x) - u^0(x) \le C.
 \end{equation}
On the other hand
\begin{eqnarray}\label{esti3}
q_j\cdot \tildex - \log l_{\xi}(y) &\ge& \frac 12 
\sum_{i \notin \{i_1,\cdots,i_m\}} l_i(q_j) \log l_i(y) \\
&&\quad- \frac 12 \sum_{i \notin \{i_1,\cdots,i_m\}} l_i(q_j) (\log |y| + \log C) + 1.\nonumber 
\end{eqnarray}
For each hyperplane $L_i = \{l_i(y) = 0\}$ we take a hyperplane $L^{\prime}_i = \{l^{\prime}_i(y) 
= 0\}$ which is close to $L_i^{\prime}$ such that
$(C(\mu)-\mathrm{apex})\cap L_i^{\prime}$ is non-empty and included in $\{l_i(y) > 0\}$
and that $(C(\mu)-\mathrm{apex})\cap L_i$ is included in $\{l_i^{\prime}(y) < 0\}$. Put
$$D_i := \{y\in C(\mu)\,|\,l_i^{\prime}(y) \le 0\}.$$
Define the following sets successively:
$$\begin{array}{l}
C_0 = C(\mu) - \bigcup_{j=1}^d D_j,\\
C_i = D_i - D_i \cap (\bigcup_{j\ne i}D_j),\\
C_{i_1i_2} = D_{i_1}\cap D_{i_2} - D_{i_1}\cap D_{i_2}\cap (\bigcup_{j\ne i_1,i_2}D_j),\\
\cdots\\
C_{i_1i_2\cdots i_k} = \bigcap_{n=1}^kD_{i_n} - \bigcap_{n=1}^kD_{i_n} \cap (\bigcup_{j\ne i_1,
\cdots,i_k}D_j),\\
\cdots\\
C_{i_1\cdots i_m} = \bigcap_{n=1}^m D_{i_n}
\end{array}$$
The union of all these sets is $C(\mu)$. 
We exclude from above the empty sets. First of all, we have on $C_0$
\begin{equation}\label{l_j}
c|y| \le l_j(y) \le C|y|
\end{equation}
for any $j$. Hence for any $q_s \in L_{j_1}\cap \cdots \cap L_{j_m}$ we have
\begin{eqnarray}\label{esti4}
&&q_s\cdot \tildex - \log l_{\xi}(y) \\
&&\ge \frac 12 \sum_{j\ne j_1,\cdots,j_m} l_j(q_s)\log l_j(y)
- \frac 12 \sum_{j\ne j_1,\cdots,j_m}l_j(q_s) (\log|y| + C) + 1\nonumber\\
&&\ge \frac 12 \sum_{j\ne j_1,\cdots,j_m} l_j(q_s)(\log |y| + c)
- \frac 12 \sum_{j\ne j_1,\cdots,j_m}l_j(q_s) (\log|y| + C) + 1\nonumber\\
&&\ge c.\nonumber
\end{eqnarray}
On $C_i$ (\ref{l_j}) holds 
for any $j \ne i$.
Take a vertex $q_s \in L_i\cap L_{j_1} \cap \cdots \cap L_{j_{m-1}}$. Then
\begin{eqnarray}\label{esti5}
&&q_s\cdot \tildex - \log l_{\xi}(y) \\
&&\ge \frac 12 \sum_{j\ne i, j_1,\cdots,j_{m-1}} l_j(q_s)\log l_j(y)
- \frac 12 \sum_{j\ne i,j_1,\cdots,j_{m-1}}l_j(q_s) (\log|y| + C) + 1\nonumber\\
&&\ge \frac 12 \sum_{j\ne i,j_1,\cdots,j_{m-1}} l_j(q_s)(\log |y| + c)
- \frac 12 \sum_{j\ne i,j_1,\cdots,j_{m-1}}l_j(q_s) (\log|y| + C) + 1\nonumber\\
&&\ge c.\nonumber
\end{eqnarray}
Continuing this way, on $C_{i_1i_2\cdots i_k}$ (\ref{l_j}) holds 
for any $j \ne i_1,\cdots, i_k$.
Take a vertex $q_s \in L_{i_1}\cap\cdots \cap L_{i_k}\cap L_{j_1} \cap \cdots \cap L_{j_{m-k}}$. Then
\begin{eqnarray}\label{esti6}
q_s\cdot \tildex - \log l_{\xi}(y) 
&\ge& \frac 12 \sum_{j\ne i_1,\cdots,i_k
j_1,\cdots,j_{m-k}} l_j(q_s)\log l_j(y) \\
&& \qquad - \frac 12 \sum_{j\ne i_1,\cdots,i_k,j_1,\cdots,j_{m-k}}l_j(q_s) (\log|y| + C) + 1
\nonumber\\
&\ge& \frac 12 \sum_{j\ne i_1,\cdots,i_k,j_1,\cdots,j_{m-k}} l_j(q_s)(\log |y| + c)\nonumber\\
&&\qquad- \frac 12 \sum_{j\ne i_1,\cdots,i_k,j_1,\cdots,j_{m-k}}l_j(q_s) (\log|y| + C) + 1
\nonumber\\
&\ge& c.\nonumber
\end{eqnarray}
It follows from (\ref{esti4}), (\ref{esti5}) and (\ref{esti6}) that
\begin{eqnarray}\label{below}
2\barv(x) - 2u^0(x) &=& 2\max_{1\le i \le \ell}  p_i\cdot x - \log l_{\xi}(y) \\
&\ge& q_s \cdot \tildex -  \log l_{\xi}(y) \ge c. \nonumber
\end{eqnarray}

Then (\ref{above}) and (\ref{below}) give the desired estimate. This completes the proof.
\end{proof}

\begin{lemma}\label{lem49}Let $X_i= - \frac i2\p\slash \p w^i$. Then we have 
\[\t_{X_i}= \frac {\p u^0}{\p x^i}.\]
\end{lemma}
\begin{proof} It is easy to show
$$ i\left(- \frac {\sqrt{-1}}2 \p\slash \p w^i\right)
\sqrt{-1}\frac{\p^2u^0}{\p x^i\p x^j}dw^i\wedge d\barw^j
= \barpartial \left(\frac{\p u^0}{\p x^i}\right).$$
Thus it is sufficient to show
$$ \int_S \frac{\p u^0}{\p x^i}\,e^h(\frac12 d\eta)^m\wedge \eta = 0.$$
But $\rho^T - \frac12 d\eta = i \p_B\barpartial_B h$ implies 
$$ e^h (\frac12 d\eta)^m = e^{-u^0}dx^1\wedge \cdots \wedge dx^m\wedge d\theta^1
\wedge \cdots \wedge d\theta^m.$$
Since $\xi$ generates isometries and $L_{\xi}\eta = 0$, we have 
$$ \int_S \frac{\p u^0}{\p x^i}\,e^h(\frac12 d\eta)^m\wedge \eta = \mathrm{const}
\int_{\bfR^m} \frac{\p u^0}{\p x^i} e^{-u^0} dx = 0$$
because $u^0$ is strictly convex and $u^0 \to \infty$ as $x \to \infty$ by 
Proposition \ref{pro46}.
\end{proof}

\begin{lemma}\label{lem50} If the vector field obtained in Proposition
\ref{pro49} is denoted by $X=\sum_{i=1}^m c_i(-\frac i2 \frac \p{\p w^i})$, Then we have
\[\int_{\Sigma }y_j e^{\sum_{i=1}^m c_i y_i} d y=0,\quad \forall j=1,2,\cdots,m.\]
\end{lemma}
\begin{proof}
We have
\[\begin{array}{rcl}
f_X(X_j)&=&\ds\vs -\int_\SS \t_{X_j} e^{\t_X} (\frac12 d\eta)^m\wedge \eta \\
&=& \ds\vs 
-\int_\SS \frac {\p u^0}{\p x^j} e^{\sum_i c_i \frac {\p u^0}{\p x^i}} \det (u^0_{ij}) dx\wedge 
d\t \wedge \eta \\
&=&-\ds \mathrm{const.} \int_{\Sigma} y_j e^{\sum_{i=1}^m c_i y_i} d y.\end{array}\]
Hence the Lemma follows from Proposition \ref{pro49}.
\end{proof}

\begin{proposition}\label{subtorus}
Let $\gamma \in \frak g^{\ast}$ be as in Proposition \ref{tangent}, and $H$ be the subtorus of
$G = T^{m+1}$ whose Lie algebra is $\frak h := \{x\,|\,(\gamma,x) = 0\}$. 
Then there is a constant $C$ such that
\begin{equation}\label{logdet}
|\log \det (u^0_{ij}) + (2m+2)u^0| \le C.
\end{equation}
\end{proposition}
\begin{proof}
We keep the same notations as before. Let $x^i$ be the real part of the  affine logarithmic 
coordinates on $Orb_{C(S)}(H^c,p) \cong (\bfC^{\ast})^m$. The natural inclusion
$Orb_{C(S)}(H^c,p) \cong (\bfC^{\ast})^m \to \mu^{-1}(\mathrm{Int}(C(\mu))$ is
induced by $j:\frak h \to \frak g$. So, we denote by $\tildex = jx$ the real part of the
affine logarithmic coordinates on 
$\mu^{-1}(\mathrm{Int}(C(\mu))\cong (\bfC^{\ast})^{m+1}$
corresponding to $x$. Let $y$ be the Legendre transform on $C(S)$ of $\tildex$.
Then the Legendre transform $v$ of $x$ on the K\"ahler manifold $Orb_{C(S)}(H^c,p)$
is given by
\begin{equation}\label{v}
v = \frac{j^{\ast}y}{l_{\xi}(y)}
\end{equation}
where $j^{\ast} : \frak g^{\ast} \to \frak h^{\ast}$. Define $\l'_i \in \frak h \cong 
\mathbb{R}^m$ by the
decomposition
\begin{equation}\label{v2}
\l_i=j\l'_i+\mu_i\xi.
\end{equation}
By our choice of $\gamma$ and $\frak h$ we have
\begin{equation}\label{v3}
\mu_i = \frac{l_i(\gamma)}{l_{\xi}(\gamma)} = \frac 1{m+1}.
\end{equation}
Moreover 
\begin{eqnarray}\label{v4}
\frac{l_i(y)}{l_\xi(y)} &=& \frac{(j\l'_i+\mu_i\xi,y)}{l_\xi(y)}\\
&=& (\l'_i,v) + \mu_i = l'_i(v)\nonumber 
\end{eqnarray}
where we have set $l'_i(v) = (\l'_i,v) + \mu_i)$. Similarly if we set
$\l'_{\infty} = \sum_{i=1}^d \l'_i$, $\mu_{\infty} = \sum_{i=1}^d \mu_i (= \frac d{m+1})$,
and $l'_{\infty}(v) = (\l'_{\infty},v) + \mu_{\infty}$ then
\begin{equation}\label{v5}
\frac{l_{\infty}(y)}{l_\xi(y)}  = l'_{\infty}(v).
\end{equation}

We see from (\ref{x^j}), (\ref{v4}) and (\ref{v5}) that 
\begin{eqnarray}\label{v6}
u^0(x) &=& \frac 12 \log l_{\xi}(y) \\
&= & -\frac12\left(\sum_{i=1}^d \mu_i\log 
\frac{l_i(y)}{l_\xi(y)}
-\mu_\infty
\log \frac{l_\infty(y)}{l_\xi(y)}
+1\right)\nonumber\\
&=& -\frac12\left(\sum_{i=1}^d \mu_i\log l'_i(v)
-\mu_\infty
\log l'_{\infty}(v)
+1\right)\nonumber
\end{eqnarray}
The symplectic potential $G_0$ on $Orb_{C(S)}(H^c,p)$ is also computed 
using (\ref{x^j}), (\ref{v4}) and (\ref{v5}) as
\begin{eqnarray}\label{v7}
G_0(v) &=& (v,x) - u^0(x)\\
&=& \frac{(j^{\ast}y,x)}{l_{\xi}(y)} - \frac 12 \log l_{\xi}(y)\nonumber\\
&=& \frac 1{l_{\xi}(y)}(y, \frac12\sum_{i=1}^d \l_i\log l_i(y)
+ \frac12 \xi(1 + \log l_{\xi}(y))
\nonumber\\
&& \qquad\qquad - \frac12\l_{\infty}\log l_{\infty}(y)) - \frac12 \log l_{\xi}(y) \nonumber\\
&=& \frac12\sum_{i=1}^d\frac {l_i(y)}{l_{\xi}(y)}\log \frac{l_i(y)}{l_{\xi}(y)}
+ \frac12 (1 + \log l_{\xi}(y))\nonumber\\
&& \qquad\qquad- \frac12\frac{l_{\infty}(y)}{l_{\xi}(y)}\log \frac{l_{\infty}(y)}{l_{\xi}(y)}
 - \frac12 \log l_{\xi}(y) \nonumber\\
&=& \frac12(\sum_{i=1}^d l'_i(v)\log l'_i(v) -  l'_{\infty}(v) \log l'_{\infty}(v) + 1)\nonumber
\end{eqnarray}
Since 
$$ \overline{\Sigma} = \bigcap_{i=1}^d \{l'_i(v) \ge 0\}$$
we have
\begin{eqnarray}\label{v8}
\det\mathrm{Hess} \,u^0(x) &=& (\det \mathrm{Hess} \,G_0(v))^{-1} \\
&=& \delta(v) \Pi_{i=1}^d l'_i(v) \nonumber
\end{eqnarray}
where $\delta$ is a strictly positive function on $\overline{\Sigma}$.
Since $\mu_i = \frac 1{m+1}$ by (\ref{v3}) we have
\begin{eqnarray*}
 |\log \det \mathrm{Hess}\, u^0 + (2m+2)u^0| &\le& 
 |\log \delta + \frac 12 \mu_{\infty}\log l'_{\infty}(v) - \frac 12| \\
&\le& C.
\end{eqnarray*}
\end{proof}

\begin{remark}\ \ One can also show Proposition 7.3 along the line of the proof of
Proposition 7.6.
\end{remark}

From now on we choose the subtorus $H$ as in Proposition \ref{subtorus}.
Then we have
\begin{equation}\label{v9}
 \det (u^0_{ij}) = \exp(-h - u^0).
\end{equation}
Now we can have our equation on $\bfR^m$. Let $\tilde \eta=\eta+2d_B^c\varphi$ be
the solution to (\ref{SR5}) where the initial metric is chosen so that the symplectic potential
on $C(S)$ is $G^{can}_{\xi}$.
Then using (\ref{SR5}) and (\ref{v9})
one can show that $u$ satisfies
\begin{equation}\label{main}
\det(u_{ij})=\exp\left(-tu -(1-t)u^0 - 
\sum_i c_i\frac{\p u}{\p x^i}\right) \quad \hbox{ on } {\mathbb R}^m.
\end{equation}
Then by the same arguments as Lemma 3.1 - 3.4 in Wang-Zhu (\cite{Wang-Zhu}) 
using our Proposition \ref{pro46} we get an estimate
$$ \sup_{Orb_{C(S)}(H^c, p)} \varphi\,\le C $$
for some constant $C > 0$ independent of $t \in [0, 1]$, or equivalently
$$ \sup_{Orb_{S}(H^c, \barp)} \varphi\,\le C.  $$
But the Reeb flow on $S$ generates isometries and we have $C^0$ estimate on
an open dense subset of $S$, which gives a $C^0$ estimate on $S$ naturally.
By the same arguments as in either
proof (i) or proof (ii) in \cite{Wang-Zhu} give an estimate 
$$ \inf_S \varphi \, \ge - C$$
for some constant $C >0$ independent of $t$. In fact we can give
all necessary modifications to the arguments of proof (ii) in \cite{Wang-Zhu},
including the arguments of Cao-Tian-Zhu (\cite{ctz}) and Mabuchi (\cite{mabuchi02}, 
\cite{mabuchi03}), which we do not re-produce here. 
This completes the proof of Lemma \ref{C^0}, and consequently
the proof of Theorem \ref{MainThm}.

\section{The invariant for K\"ahler cone manifolds}

In this section, we will reformulate the invariant $f = 2\pi f_1$
of Sasaki manifolds of positive basic first Chern class
as an invariant for K\"ahler cone manifolds,
and then relate the volume function of Sasaki manifolds with 
the invariant $f$. This relation was pointed out in section 4.2 of \cite{MSY2}
in the case of quasi-regular Sasaki manifolds.
We wish to relate the invariant  further to the existence problem of a Sasaki-Einstein metric.

In the previous sections we used the same notation $\xi$ for the
Reeb field $J\frac \p{\p r}$ and the vector field $Jr\frac{\p}{\p r}$ on $C(S)$,
but we distinguish them by denoting the vector field on $C(S)$ as
$$ \widetilde{\xi} = Jr\frac{\p}{\p r}.$$

\begin{defn}\label{11}
Let $S$ be a compact $(2m+1)$-dimensional manifold and 
a complex structure $J$ on $C(S):=\mathbb{R}_+\times S$.
We call a K\"ahler metric $\bar{g}$ on $(C(S),J)$ a K\"ahler cone metric
if there exist a smooth function $r:C(S)\to \mathbb{R}_+$, a smooth map
$p:C(S)\to S$ and a Riemannian metric $g$ on $S$ such that
$(r,p):C(S)\to \mathbb{R}_+\times S$ is a diffeomorphism and that 
$(r,p)^*(ds^2+s^2g)=\bar{g}$. (Then, of course, $g$ is a Sasaki metric on $S$.)
\end{defn}

When a K\"ahler cone metric $\bar{g}$ on $(C(S),J)$ is given, 
we define the vector field $\widetilde{\xi} $ and the $1$-form $\eta$
on $C(S)$ by
$$\widetilde{\xi} =rJ\frac{\partial}{\partial r},\ \eta=
\frac{1}{r^2}\bar{g}(\widetilde{\xi} ,\cdot)=2d^c\log r$$
where we use the notation $d^c = \frac{ i}{2} (\barpartial - \p)$. 
Then we see that $\widetilde{\xi} $ is a holomorphic and Killing vector field.
Moreover $\widetilde{\xi} $ lies in the center of the Lie algebra of 
the group of isometries $\isom(C(S),\bar{g})$.
The restrictions of $\widetilde{\xi} $ and $\eta$ to
$\{r=1\}\subset C(S)$, where $r$ is the smooth function on $C(S)$ associated
with the K\"ahler cone metric $\bar{g}$, are the Reeb vector field and
the contact $1$-form on the Sasaki manifold $\{r=1\}\simeq (S,g)$.
Moreover we see that the K\"ahler form $\omega$ of $\bar{g}$ has 
the K\"ahler potential 
$\frac12r^2$;
\begin{equation}
\omega=\frac12d(r^2\eta)=\frac{i}{2}\partial\bar{\partial}r^2.
\end{equation}

Let $S$ be a compact $(2m+1)$-dimensional manifold and 
$J$ be a given complex structure on $C(S)$.
Suppose that the canonical bundle $K_{C(S)}$ of $C(S)$ is trivial.
Then we want to know whether a Ricci-flat K\"ahler cone metric exists on
$(C(S),J)$. In what follows we will reformulate the invariant $f$ obtained in 
Theorem \ref{Futaki-inv} as an 
obstruction to the existence of a  Ricci-flat K\"ahler cone metric on $C(S)$.

We fix a maximal torus $T^n\subset \aut(C(S),J)$. (Later we will consider toric Sasaki manfolds,
and then $n = m+1$. But for the moment we do not assume $S$ is toric.)
Let $\KCM(C(S),J)$ denote
the space of K\"ahler cone metrics on $(C(S),J)$
such that the maximal torus $T^n$ is contained in the group of isometries
and $\widetilde{\xi} \in \mathfrak g$, where $\mathfrak g$ is the Lie algebra of $T^n$.
For each metric in $\KCM(C(S),J)$, there is an associated moment map
$$\mu:C(S)\to \mathfrak g^*,\ \ \langle \mu,X\rangle= r^2\eta(X)$$
where we identify $X \in \mathfrak g$ with the corresponding vector field on $C(S)$
(recall our convention of the moment map in section 6 where we deleted $\frac 12$).
The image of $\mu$ is a convex rational polyhedral cone $C(\mu)
\subset {\mathfrak g}^*$. 
Moreover these cones are all isomorphic for all K\"ahler cone metrics in $\KCM(C(S),J)$.
Note that, since the K\"ahler form
of any metric in $\KCM(C(S),J)$ is exact, this is a deformation of 
K\"ahler metrics with the same K\"ahler class.

When we investigate the existence of Ricci-flat K\"ahler cone metrics,
we may restrict the deformation space of K\"ahler cone metrics to
$$\KCM_c(C(S),J):=\{\bar{g}\in \KCM(C(S),J)\,|\,
\rho(\bar{g})=i\partial \bar{\partial}\tilde{f},\ 
r\frac{\partial}{\partial r}\tilde{f}=\widetilde{\xi} \tilde{f}=0\}.$$
Of course a Ricci-flat K\"ahler cone metric belongs to $\KCM_c(C(S),J)$.
In the Sasakian point of view, $\bar{g}\in \KCM_c(C(S),J)$ means
that $[\rho^T]_{B,\xi}=2(m+1)[\omega^T]_{B,\xi}$, where $[\rho^T]_{B,\xi}$ and 
$[\omega^T]_{B,\xi}$
are respectively the basic cohomology classes of the transverse Ricci form and 
the transverse K\"ahler form of 
the Sasaki metric $g$ on $S$ induced from $\bar{g}$, the Reeb field $\xi$ being with
respect to $g$.
This is because $\rho(\bar{g})=\rho^T-2(m+1)\omega^T$ 
and $r\frac{\partial}{\partial r}\tilde{f}=\widetilde{\xi} \tilde{f}=0$. 
Proposition \ref{tangent} can be restated as follows.

\begin{proposition}\label{tangent2}
$\KCM_c(C(S),J)$ is exactly the set of
all K\"ahler cone metrics such that 
the associated Reeb fields satisfy (\ref{F4}), or equivalently (\ref{F9}).
\end{proposition}

We also consider the subset 
$$\KCM_c(C(S),J,\widetilde{\xi} _0):=
\{\bar{g}\in \KCM_c(C(S),J)\,|\,\widetilde{\xi} =\widetilde{\xi} _0\}$$
for each fixed $\widetilde{\xi} _0$. This subset 
corresponds to the set of all the transverse K\"ahler deformations of Sasaki metrics
with the fixed Reeb vector field $\widetilde{\xi} _0$ and varying isometric inclusion $S\subset C(S)$.
This fact can be checked as follows.
Let $\bar{g}$ and $\bar{g'}$ be K\"ahler cone metrics on
$(C(S),J)$ with $\widetilde{\xi} =\widetilde{\xi} '  = \widetilde{\xi} _0$. 
Then by rotating by $J$ we get
$r\frac{\partial}{\partial r}=r'\frac{\partial}{\partial r'}$. 
Hence there exists a smooth function $\varphi$
on $C(S)$ such that $r'=r\exp(\varphi)$ and $\frac{\partial}{\partial r}\varphi=\widetilde{\xi} \varphi=0$.
On the other hand, $(S,g')$ is isometrically identified with $\{r'=1\}
\subset C(S)$
by the embedding $i:S\hookrightarrow C(S)$, 
$i(x)=(\exp(-\varphi)(x),x)$. Thus, on $\{r'=1\}$,
$i_*(\xi)=0\oplus \xi=\widetilde{\xi} _
{\{r'=1\}}=\widetilde{\xi} '_
{\{r'=1\}}$.
Since $i_*$ is injective, this suggests that $\xi$ is the Reeb vector field
of the Sasaki manifold $(S,g')$.
The one form $\eta$ on $C(S)$ 
is deformed to 
$$\eta'=2d^c(\log r+\varphi)=\eta+2d^c\varphi$$
when we change the metric from $\bar{g}$ to $\bar{g'}$.
Hence 
$$\omega^T(g')=\omega^T(g)+d_Bd_B^c\varphi,$$
where $\omega^T(g),\omega^T(g')$ are the transverse K\"ahler forms 
with respect to Sasaki metrics $g,g'$ on $S$ respectively.

Consider the volume functional
$$\tilde{V}:\KCM(C(S),J)\to \mathbb{R},\ \ \tilde{V}(\bar{g}):=\vol(S,g),$$
where $g$ is the Sasaki metric on $S$ induced from $\bar{g}$.
Proposition \ref{19} and \ref{15} below, which are the first and second variation formulae for the volume functional,
were proved in Appendix C of \cite{MSY2},
but we give slightly more comprehensive proofs in this paper for the reader's
convenience.

\begin{prop}[\cite{MSY2}, Appendix C]\label{19}
Let $S$ be a compact Sasaki manifold.
Let $\{\bar{g}(t)\}_{-\varepsilon<t<\varepsilon}$ be a $1$-parameter family
in $\KCM(C(S),J)$ with $\bar{g}(0)=\bar{g}$ and $Y=d\widetilde{\xi} /dt_{|t=0}$.
Then 
\begin{equation}\label{20}
\frac{d}{dt}_{|t=0}\tilde{V}(\bar{g}(t))=-4(m+1)
\int_S\eta(Y)dvol_g.
\end{equation}
\end{prop}
\begin{proof}
Let $S(t)\subset C(S)$ denote the subset $\{r(t)=1\}$ for each $t$,
where $r(t):C(S)\to \mathbb{R}_+$ is the smooth function associated
with the cone metric $\bar{g}(t)$, see Definition \ref{11}. 
Then $(S(t),\bar{g}(t)
_{|S(t)})$ is isometric to $(S,g(t))$.
Since $(d\eta(t))^m\wedge \eta(t)$ is closed for each
$t$,
we have
\begin{align*}
\tilde{V}(\bar{g}(t)) &= \frac{1}{m!}\int_{S(t)}(\frac12 d\eta(t))^m\wedge 
\eta(t)=\frac{1}{m!}\int_{S(0)}(\frac12 d\eta(t))^m\wedge 
\eta(t)\\
&= \frac{1}{m!} \int_{S(0)}\left(\frac12 d\eta+ t dd^c\phi+O(t^2)\right)^m
\wedge \left(\eta+2t d^c\phi+O(t^2)\right)
\end{align*}
Therefore, by Lemma \ref{17}, the first variation is
\begin{align*}
\frac{d}{dt}_{|t=0}\tilde{V}(\bar{g}(t)) &=
\frac{1}{m!}\int_{S(0)}\bigg(m(\frac 12 d\eta)^{m-1}\wedge 
dd^c\phi\wedge\eta+(\frac12 d\eta)^m\wedge 2d^c\phi\bigg)\\
&= \frac{2m+2}{m!}\int_{S(0)}d^c\phi\wedge (\frac12 d\eta)^m
=-4(m+1)\int_S\eta(Y)dvol_g.
\end{align*}
Here the second equality holds since
$$d((d\eta)^{m-1}\wedge d^c\phi\wedge \eta)=
(d\eta)^{m-1}\wedge 
dd^c\phi\wedge\eta-(d\eta)^m\wedge d^c\phi.
$$
\end{proof}
\begin{prop}[\cite{MSY2}, Appendix C]\label{15}
Let $\{\bar{g}(t)\}_{-\varepsilon<t<\varepsilon}$ be a $1$-parameter family
in $\KCM_c(C(S),J)$ with $\bar{g}(0)=\bar{g}$ and $Y=d\widetilde{\xi} /dt_{|t=0}$.
Then 
\begin{equation}
\frac{d}{dt}_{|t=0}\left(\int_S\eta(t)(X)dvol_{g(t)}\right)=-(2m+4)
\int_S\eta(X)\eta(Y)dvol_g
\end{equation}
for each $X\in \mathfrak g$, where $g(t)$ and $\eta(t)$ are
the Sasaki metric and the contact $1$-form on $S$ induced from
$\bar{g}(t).$
\end{prop}

\begin{proof}
Let the notations be as in the proof of the previous proposition.
For each sufficiently small $t$, 
\begin{equation}\label{e14}
\widetilde{\xi} (t)=\widetilde{\xi} +tY+O(t^2)
\end{equation}
and
\begin{equation}\label{e15}
r^2(t)=r^2(1+2t\varphi+O(t^2)).
\end{equation}
Then we have
\begin{equation}\label{e16}
\eta(t)=\eta+2t d^c\varphi+O(t^2).
\end{equation}
Multiplying \eqref{e14} by $J$, we have
\begin{equation}\label{e17}
r(t)\frac{\partial}{\partial r(t)}=r\frac{\partial}{\partial r}-tJY+O(t^2).
\end{equation}
Expanding $\mathcal{L}_{r(t)\frac{\partial}{\partial r(t)}}r^2(t)
=2r^2(t)$ to first order in $t$ gives
\begin{equation}\label{e18}
2\mathcal{L}_{r\frac{\partial}{\partial r}}\varphi=-2\eta(Y).
\end{equation}
Since $\mathcal{L}_{r(t)\frac{\partial}{\partial r(t)}}X=0$ for each $X\in \mathfrak g$
and $\mathcal{L}_{r(t)\frac{\partial}{\partial r(t)}}\eta(t)=0$,
\begin{align*}
d\left\{\eta(t)(X)\left(d\eta(t)\right)^m\wedge
\eta(t)\right\} &= 
d\left(\eta(t)(X)\right)\wedge\left(d\eta(t)\right)^m\wedge
\eta(t)\\
&= \frac{\partial (\eta(t)(X))}{\partial r(t)}dr(t)
\wedge\left(d\eta(t)\right)^m\wedge
\eta(t)=0
\end{align*}
on $C(S)$. Hence
\begin{align*}
\int_S\eta(t)(X)dvol_{g(t)} &= \frac{1}{m!}\int_{S(t)}\eta(t)(X)
\left(\frac12d\eta(t)\right)^m\wedge
\eta(t)\\
&= \frac{1}{m!}\int_{S(0)}\eta(t)(X)
\left(\frac12d\eta(t)\right)^m\wedge\eta(t)\\
&= \frac{1}{m!}\int_{S(0)}\bigg\{
\left(\eta+2t d^c\varphi+O(t^2)\right)
(X)\left(\frac12 d\eta+ t dd^c\varphi+O(t^2)\right)^m\\
&\ \ \ \ \ \ \ \ \wedge
\left(\eta+ 2t d^c\varphi+O(t^2)\right)\bigg\}.
\end{align*}
Thus
\begin{align*}
\frac{d}{dt}_{|t=0}\left(\int_S\eta(t)(X)dvol_{g(t)}\right) &=
\frac{1}{m!}\int_{S(0)}\bigg\{2d^c\varphi(X)(\frac12 d\eta)^m\wedge
\eta\\
&\ \ \ \ \ +m\eta(X)(\frac12 d\eta)^{m-1}\wedge dd^c\varphi
\wedge\eta\\
&\ \ \ \ \ +\eta(X)(\frac12 d\eta)^m\wedge d^c\varphi\bigg\}\\
&=\frac{2m+4}{m!}\int_{S(0)}\eta(X)(\frac12 d\eta)^m\wedge 2 d^c\varphi\\
&= -(2m+4)
\int_S\eta(X)\eta(Y)dvol_g.
\end{align*}
Here the second and third equalities are given by the following lemmas.
\end{proof}

\begin{lem}\label{16}
\begin{eqnarray*}
&&\int_{S(0)}\left\{d^c\varphi(X)(\frac12 d\eta)^m\wedge
\eta+\frac m2 \eta(X)(\frac12 d\eta)^{m-1}\wedge dd^c\varphi
\wedge\eta\right\}\\
&& =(m+1)\int_{S(0)}\eta(X)(\frac12 d\eta)^m\wedge d^c\varphi.
\end{eqnarray*}
\end{lem}

\begin{proof}
Since $X$ is tangent to $S(0)$ and $\mathcal{L}_X\eta=0$,
\begin{eqnarray*}
0 &=& \int_{S(0)}\iota(X)((\frac12 d\eta)^m\wedge \eta
\wedge d^c\varphi)\\
&=& \int_{S(0)}(m(\frac12 \iota(X)d\eta)\wedge (\frac12 d\eta)^{m-1}
\wedge \eta\wedge d^c\varphi\\
&& \quad +\eta(X)(\frac12 d\eta)^m
\wedge d^c\varphi
-d^c\varphi(X)(\frac12 d\eta)^m\wedge \eta)\\
&=& \int_{S(0)}(-\frac m2 d(\eta(X))\wedge (\frac12d\eta)^{m-1}
\wedge \eta\wedge d^c\varphi\\
&&\quad+\eta(X)(\frac12 d\eta)^m
\wedge d^c\varphi-d^c\varphi(X)(\frac12 d\eta)^m\wedge \eta).
\end{eqnarray*}
On the other hand 
\begin{eqnarray*}
&&\frac12 d\big(\eta(X)(\frac12 d\eta)^{m-1}
\wedge \eta\wedge d^c\varphi\big) \\
&&\qquad= \frac12 d(\eta(X))\wedge (\frac12 d\eta)^{m-1}
\wedge \eta\wedge d^c\varphi+\eta(X)(\frac12 d\eta)^m\wedge
d^c\varphi\\
&&\qquad\qquad - \frac12 \eta(X)(\frac12 d\eta)^{m-1}\wedge dd^c\varphi
\wedge\eta.
\end{eqnarray*}
Combining these equations, we get the lemma.
\end{proof}

\begin{lem}\label{17}
On $S(0)$,
\begin{equation}\label{e110}
d^c\varphi\wedge (d\eta)^m=-2\eta(Y)\eta\wedge
(d\eta)^m.
\end{equation}
\end{lem}

\begin{proof}
On $S(0)$,
$$d^c\varphi\wedge (d\eta)^m=(\iota(\widetilde{\xi} )d^c\varphi)\eta
\wedge (d\eta)^m.$$
On the other hand, by \eqref{e18},
$$\iota(\widetilde{\xi} )d^c\varphi=-2\iota(J\widetilde{\xi} )d\varphi=2\iota\left(
r\frac{\partial}
{\partial r}\right)d\varphi=-2\eta(Y).$$
Combining these equations, we get \eqref{e110}.
\end{proof}

Let $\bar{g}\in \KCM_c(C(S),J)$, and denote by
$\mathfrak{h}(C(S),J,\bar{g})$ the space of {\it normalized} holomorphic Hamiltonian 
vector fields on $(C(S),J,\bar{g})$ in the following sense.
We call a complex vector field $\tildeX$ on $C(S)$ Hamiltonian holomorphic 
if $\tildeX$ is a holomorphic vector field and $\tildeX_\mathbb{R}=(\tildeX+\overline{\tildeX})/2$
is Killing. If $\tildeX$ is a holomorphic Hamiltonian vector field then
$$X=\widetilde{X}-i\eta(\widetilde{X})r\frac{\partial}{\partial r}
=\widetilde{X}-i\eta(X)r\frac{\partial}{\partial r}$$
defines a Hamiltonian holomorphic vector field on $S$ in the sense of Definition 4.4.
With this remark, we say that a Hamiltonian holomorphic 
vector field $\tildeX$ on $C(S)$ is normalized if $X$ is normalized in the sense that
$u_X := \frac i2 \eta(X)$ satisfies (\ref{normal}).

Then, using the above relation between $\tildeX$ and $X$, we define a linear function
on $\mathfrak{h}(C(S),J,\bar{g})$ by
$$\tildeX\mapsto F(\tildeX):= 2\pi i\int_S X h\,dvol_g,$$
where $g$ is the Sasaki metric on $S$ induced from the cone metric $\bar{g}$
and $h$ is the basic function on $S$ such that $\widetilde{f}=p^*h$.

\begin{prop}\label{12}
The linear map $F$ defined above coincides with 
the restriction to $\mathfrak{h}(C(S),J,\bar{g})$ of  $2\pi im!f = m! if_1$
on $(S,g)$.
\end{prop}

\begin{proof} This follows from (\ref{f}).
\end{proof}

\begin{prop}\label{13}
$F$ is independent of $\bar{g} \in
\KCM_c(C(S),J,\widetilde{\xi} _0)$
for each fixed $\widetilde{\xi} _0$. 
\end{prop}
\begin{proof} This follows from Proposition \ref{12} and Theorem \ref{Futaki-inv}.
But we will give an alternate proof below.
\end{proof}

Notice that $\mathfrak{h}(C(S),J,\bar{g})$ may vary as $\xi$ varies as
the elements of $\mathfrak{h}(C(S),J,\bar{g})$ have to commute with $\xi$.
Since this linear function $F$ is a character of the Lie algebra
$\mathfrak{h}(C(S),J,\bar{g})$, $F$ is only nontrivial on the center, and
$\xi$ and the center are included in the maximal torus $\mathfrak g$.
So we restrict $F$ to $\mathfrak g$, but consider it for all $\bar{g} \in \KCM_c(C(S),J)$:
$$ F : \KCM_c(C(S),J) \times \mathfrak g \to \bfR.$$

\begin{proof}[Alternative Proof of Proposition \ref{13}]
Since $X$ is a normalized Hamiltonian holomorphic vector field, $\eta(X)=-  2 i u_X $
satisfies 

\begin{eqnarray}\label{e12}
2(m+1)\pi \int_S \eta(X)\, dvol_g &=&
 \pi \int_S(- \Delta_B \eta(X) 
+2(m+1)\eta(X))\,dvol_g  \\
&=& -  2\pi i\int_S \nabla^i u_X \nabla_ih\,dvol_g\nonumber\\
&=& 2\pi i\int_S Xh\,\,dvol_g\nonumber\\
&=& F(\tildeX).\nonumber
\end{eqnarray}
Therefore, 
the invariance of $F$ on $\KCM_c(C(S),J,\xi_0)$
follows from the following Proposition \ref{15} by putting $Y = 0$.
\end{proof}

Of course if $\bar{g}\in \KCM_c(C(S),J)$ is 
Ricci-flat, then $F$ vanishes on $\KCM_c(C(S),J, \widetilde{\xi})$
for the corresponding $\widetilde{\xi}$. Therefore the nonvanishing of $F$ 
on $\KCM_c(C(S),J)$ obstructs
the existence of a Ricci-flat K\"ahler cone metric. 

Now let $(S,g_0)$ be a $(2m+1)$-dimensional compact toric Sasaki manifold
(see Definition \ref{def42}).
Then 
the metric cone $(C(S),J,\bar{g}_0)$ is a toric K\"ahler cone.
Here a K\"ahler metric being toric 
means that the real torus
$T^{m+1}$ acts holomorphically and effectively on $C(S)$
preserving the K\"ahler form.
Note here that if we fix a maximal torus $T^{m+1}$ of $Aut(C(S),J)$, 
then the metrics in $\KCM(C(S),J)$ are all toric K\"ahler
since we defined $\KCM(C(S),J)$ to be the set of all K\"ahler cone metric invariant
under the maximal torus of $Aut(C(S),J)$.
We will see that there is a unique vector field $\xi_c$ on
$S$ such that there are Sasaki metrics on $S$ such that
the Reeb vector field is $\xi_c$ and that the invariant $f$ 
vanishes identically. To find such $\xi_c$, we need to use the relationship
between 
the invariant $f$ and the volume functional of Sasaki manifolds
given by the first variation formula, Proposition \ref{19}.
In \cite{MSY1} and
\cite{MSY2}, Martelli, Sparks and Yau came up with this idea.

When we fix an angular coordinates $\phi_i\sim \phi_i+2\pi$ on $T^{m+1}$,
we can identify $\mathfrak g$ with $\mathbb{R}^{m+1}$ by identifying $\sum_i X^i\partial
/\partial \phi_i$ with $(X^1,\cdots,X^{m+1})$, 
and $\mathfrak g^*$ is also identified with $\mathbb{R}^{m+1}$.
Then the cone $C(\mu)$ can be represented as
$$C(\mu)=\{y\in \mathbb{R}^{m+1}|\langle \lambda_j,y\rangle\ge 0,
j=1,\cdots,d\},$$
where $\lambda_j$ are the inward pointing normal vectors of the $d$ facets of 
the cone $C(\mu)$. Since $C(\mu)$ is a rational cone,
we can normalize the vectors $\lambda_j$ to be primitive elements of $\mathbb{Z}
^{m+1}$. 

As is shown in section 2 of \cite{MSY1}, the image $\mathrm{Im}(Reeb)$ of the map
$$Reeb:\KCM(C(S),J)\to \mathfrak g \simeq \mathbb{R}^{m+1},\ \ 
\bar{g}\mapsto \widetilde{\xi} =rJ\frac{\partial}{\partial r}$$
is $C(\mu)^{\ast}_0$, the interior of the dual cone $C(\mu)^{\ast}$ of 
$C(\mu)$, see (\ref{dualcone}).
\begin{prop}\label{18}
Let $S$ be a compact toric Sasaki manifold.
The volume functional
$$\wt{V}:\KCM(C(S),J)\to \mathbb{R},\ \ \wt{V}(\bar{g}):=\vol(S,g),$$
where $g$ is the Sasaki metric on $S$ induced from $\bar{g}$,
descends to a well-defined function $V$ on $C(\mu)^{\ast}_0$
and thus depends
only on Reeb fields.
\end{prop}
This proposition is a corollary to Proposition \ref{19}:
One just puts $Y = 0$.

Combining Proposition \ref{18} and the second variation formula, 
Proposition \ref{15},
the volume function  $V:C(\mu)^{\ast}_0\to \mathbb{R}$ is a non-negative
continuous strictly convex function.
In fact, $V(x)$ can be described by the Euclidean volume of the polytope
$\Delta_x$, which depends only on $x\in C(\mu)^{\ast}_0$, in 
$\mathfrak g^*\simeq \mathbb{R}^{m+1}$;
$$\Delta_x=\{y\in C(\mu)\ |\ 2\langle x,y\rangle\le 1\},$$
$$V(x)=2(m+1)(2\pi)^{m+1}\vol(\Delta_x),$$
where $\vol(\Delta_x)$ is the Euclidean volume of the polytope
$\Delta_x$, see section 2 of \cite{MSY1}.
Hence we see that $V$ diverges to infinity
when $x\to \partial C(\mu)^{\ast}$.

By Proposition \ref{tangent}
\begin{eqnarray*}
&&\{\wt\xi \in C(S)^{\ast}_0 \,|\, \wt\xi \hbox {\ is the Reeb field for some\ } \overline{g} \in \KCM_c(C(S),J)\} \\
&&\qquad\qquad= \{\wt\xi \in C(S)^{\ast}_0 \,|\, (\gamma,\wt\xi) = -(m+1)\}
\end{eqnarray*}
which is a relatively compact set in $C(S)^{\ast} $. 
Therefore the restriction 
$V:C(\mu)^{\ast}_0\cap \{x\in \mathbb{R}^{m+1}|(\gamma,\wt\xi)=-(m+1)\}
\to \mathbb{R}$ has unique minimum point $x_c$.

\begin{prop}\label{110}
Let $(S,g_0)$ be a $(2m+1)$-dimensional compact Sasaki manifold such that
the metric cone $(C(S),J,\bar{g}_0)$ is a toric K\"ahler cone 
$\bar{g}_0\in \KCM_c(C(S),J)$. Suppose that $\widetilde{\xi} =x_c
\in C(\mu)^{\ast}_0$. Then the invariant $f$ for
$(S,g_0)$ vanishes.
\end{prop}

\begin{proof}
Let $\bar{g}\in \KCM_c(C(S),J)$ with $\widetilde{\xi} \in C(\mu)^{\ast}_0\cap \{x\in 
\mathfrak g^{\ast}\cong\mathbb{R}^{m+1}|(\gamma,x)= -(m+1)\}$.
By \eqref{e12} and the first variation formula (Proposition \ref{19})
it suffices to prove that 
$$\{X\in \mathfrak g\,|\,(\gamma,X) = 0\}=\{X\in \mathfrak g\, |
\, \Delta_B^h u_X=2(m+1)u_X\},$$
where $u_X = \frac i2 \eta(X_c)$, $X_c$ is the Hamiltonian holomorphic vector field 
$$X_c=X-i\Phi X$$
on $S$ and $h$ is a basic function which satisfies $\rho(\bar{g})=
i\partial \bar{\partial}(p^{\ast}h)$. (Note that we use the notation  $X_c$ for $X$ in Definition 4.4.)

Then by Proposition \ref{tangent}
$$\{X\in \mathfrak g\ |\ (\gamma,X)=0\}=\{X\in \mathfrak g\, |\, 
\mathcal{L}_{\tildeX_c}\Omega=0\},$$
where $\tildeX_c = X_c + i\eta(X_c)r\frac {\p}{\p r}$ and 
$\Omega$ is the non-vanishing multi-valued holomorphic $(m+1,0)$-form satisying (\ref{e112}). 
Let $X\in \{X\in \mathfrak g\ |\ \mathcal{L}_{\tildeX_c}\Omega=0\}$
be given. Then the Lie derivative of \eqref{e112} by $\widetilde{X}_c$ gives along $\{r=1\} \cong S$
\begin{equation}\label{e113}
\begin{split}
0 &=2\tildeX_ch-\Delta_{C(S)}(\frac i2 r^2\eta(\tildeX_c))\\
&= 2X_ch+\left(\frac{\partial^2}{\partial r^2} - \frac{1}{r^2}\Delta_S 
+ \frac{2m+1}{r}
\frac{\partial}{\partial r}\right)(\frac i2 r^2\eta(X_c))\\
&= 2X_ch+4(m+1)u_X-\Delta_Su_X\\
&= -2\Delta_B^hu_X+4(m+1)u_X,
\end{split}
\end{equation}
where $\Delta_{C(S)}$ and $\Delta_S$ are the positive real Laplacians of $(C(S),\bar{g})$
and $(S,g)$ respectively and we have put $u_X = \frac i2 \eta(X_c)$. Hence $\{X\in \mathfrak g\ |\ (\gamma,X)=0\}\subset\{X\in \mathfrak g\ |
\ \Delta_B^hu_X=2(m+1)u_X\}$.
Therefore $\{X\in \mathfrak g\ |\ (\gamma,X)=0\}=\{X\in \mathfrak g\ |
\ \Delta_B^h u_X=2(m+1) u_X\}$, since
these are hyperplanes in $\mathfrak g$.
Therefore $Y$ in Proposition \ref{19} can be taken as normalized 
Hamiltonian holomorphic vector fields.
\end{proof}

\begin{remark}
In the quasi-regular case, the relationship between the first variation of the 
volume function $V$ and the invariant in \cite{F0} of the orbit space of the flow of 
the Reeb vector field, which is in general K\"ahler orbifold, was proved
in section 4.2 of \cite{MSY2}, using ``Killing spinor" on $S$. We have proved this relationship
without using Killing spinors.
\end{remark}

\begin{rem}\label{111}
Note that there always exist toric Sasaki metrics such that $\widetilde{\xi} =x_c
\in C(\mu)^{\ast}_0$. In fact, we can construct symplectic potential
of such metrics concretely, see section 2 of \cite{MSY1}.
\end{rem}

\begin{proof}[Proof of Theorem \ref{thm2}]
By Proposition \ref{110} there exists $\widetilde{\xi}$ such that the corresponding
K\"ahler cone metric with vanishing invariant $f$. By Theorem \ref{MainThm} there
exists a transverse K\"ahler-Einstein metric satisfying
$\rho^T = (2m+2) \omega^T$.  This metric is a Sasaki-Einstein metric.
\end{proof}

Let $(M,\underline{g},J)$ be a real $2m$-dimensional compact K\"ahler manifold 
such that $[\rho]=2(m+1)[\omega]\in H^2(M;\mathbb{R})$, 
where $\omega$ and $\rho$ are 
the K\"ahler form and the Ricci form respectively.
Let $N$ be the maximal integer such that $c_1(M)/N$ is an integral cohomology 
class and $\pi:S\to M$ the principal $S^1$-bundle with the first Chern class
$c_1(S)=c_1(M)/N$. Then it is well-known that $S$ is simply connected and
there is a regular Sasaki metric $g$ on $S$ such that the projection $\pi$ 
is a Riemannian submersion. Moreover, this regular Sasaki metric $g$ is 
Einstein if $\underline{g}$ is Einstein. However, in contrast, $g$ is not
Einstein if $\underline{g}$ is not Einstein.

\begin{proof}[Proof of Corollary \ref{cor2}]
Let $M$ be the blow up of $\mathbb{C}P^2$ at $2$ generic points.
Then $M$ is a toric manifold which does not admit K\"ahler-Einsten 
metric by Matsushima's theorem. Hence any regular Sasaki metric on $S$
associated with K\"ahler metric on $M$ is not Einstein.

However, there is a toric Sasaki-Einstein metric on $S$ by Theorem \ref{thm2};
in this case, the inward pointing normal vectors of the facets of the 
moment cone $C(\mu)\subset \mathbb{R}^3$ of $C(S)$ are
$$v_1=(1,0,0),\ v_2=(1,0,1),\ v_3=(1,1,2),\ v_4=(1,2,1),\ v_5=(1,1,0).$$
By the calculation in section 3 of \cite{MSY1}, the Reeb vector of the toric Sasaki
Einstein metric is given by
$$x_c=\left(3,\frac{9}{16}(-1+\sqrt{33}),\frac{9}{16}(-1+\sqrt{33})\right).$$
Clearly, since $x_c$ is not a rational point, this is the Reeb vector of an
irregular Sasaki metric.
\end{proof}

\section{Appendix}

This proof of Theorem \ref{Futaki-inv} is based on the following two lemmas. 

\begin{lem}\label{lemma24} If $\a$ is a basic $(2m-1)$-form, then
\[\int_\SS d_B\a\wedge \eta=0. \]\end{lem}
This lemma can be seen as a special case of the following

\begin{lem}\label{lemma25} If $\a$ and $\b$ are basic forms with $\deg \a+
\deg \b=2m-1$, then
\[\int_\SS d_B\a\wedge\eta=(-1)^{\deg \a}\int_\SS \a\wedge d_B \b\wedge \eta.\]
\end{lem}

The proof of these lemmas are easy exercises. These lemmas
show that if we have a result for a compact K\"ahler manifold which can be proved only using
the Stokes theorem, including the integration by parts, then the result holds true for compact
Sasaki manifolds. The proof can be given only by adding "$\wedge \eta$" at each line
of the proof.

\end{document}